\documentclass[12pt]{article}
\usepackage{graphicx}
\usepackage{amsmath}
\usepackage{graphicx}
\usepackage{amsfonts}
\usepackage{amssymb}
\usepackage{amsmath, amssymb ,amsthm, amsfonts, amsgen}
\DeclareGraphicsExtensions{.eps,.bmp,.jpg,.pdf,.mps,.png,.gif}
\numberwithin{equation}{section} 
  \newtheorem{Theorem}{Theorem}[section]

\newtheorem{Proposition}[Theorem]{Proposition}
\newtheorem{Corollary}[Theorem]{Corollary}

\newtheorem{Remark}[Theorem]{Remark}

\newcommand\AMSname{AMS subject classifications}

\begin{document}\title{ Asymptotic Analysis of the Eigenvalues
of a Laplacian Problem  in a Thin Multidomain}
\author{Antonio Gaudiello\footnote{DAEIMI,
Universit\`a degli Studi di Cassino, via G. Di Biasio 43, 03043
Cassino (FR), Italia. e-mail: gaudiell@unina.it}  $\,$and Ali
Sili\footnote{D\'epartement de Math\'ematiques, Universit\'e du
Sud  Toulon-Var, BP 20132, 83957 La Garde cedex, France, \& LATP,
UMR 6632, Universit\'e de Provence, 39 rue F. Joliot-Curie, 13453
Marseille cedex 13, France. e-mail: sili@univ-tln.fr}}\date{ }
\maketitle

\begin{abstract}
We consider a thin multidomain of $\mathbb R^N$, $N\geq 2$,
consisting   of two vertical cylinders, one placed upon the other:
the first one with given height  and small cross section, the
second one with small thickness and given cross section. In this
multidomain we study the asymptotic behavior, when the volumes of
the two cylinders vanish, of a Laplacian eigenvalue problem and of
a $L^2$-Hilbert orthonormal basis of eigenvectors. We derive the
limit eigenvalue problem (which is well posed in the union of the
limit domains, with respective dimension $1$ and $N-1$) and the
limit basis.  We discuss the limit models and we precise how these
limits depend on the dimension $N$ and on limit $q$ of the ratio
between the volumes of the two cylinders.


\par
\noindent2000 \AMSname: 49R50, 35P20, 74K05, 74K30, 74K35.
\end{abstract}

\section{Introduction and main results}\label{introduction}

For every  $n \in \mathbb N$, let $\Omega_n\subset \mathbb R^N$,
$N\geq 2$, be a thin multidomain consisting of two vertical
cylinders, one placed upon the other: the first one with constant
height $1$ and small cross section $r_n\omega$, the second one
with small thickness $h_n$ and constant cross section $\omega$,
where $\omega$ is a  bounded open connected set of $\mathbb
R^{N-1}$ containing the origin $0'$ of $\mathbb R^{N-1}$ and with
smooth boundary, $r_n$ and $h_n$ are two small parameters
converging to zero. Precisely,
$$\Omega_n=\left(r_n\omega\times [0,1[\right)\bigcup\left(\omega\times
]-h_n,0[\right)$$ (for instance, see Fig. 1 for $N=2$ and Fig. 2
for $N=3$).
\begin{figure}[h] \centering
\includegraphics[width=7cm]{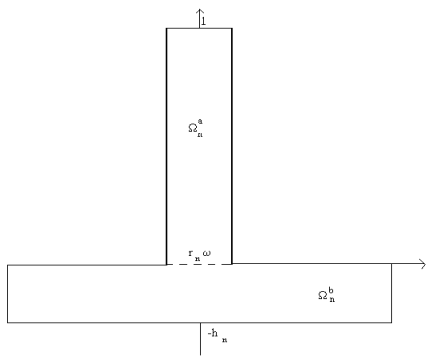}
\caption{the thin multidomain for $N=2$.} \label{domain1}
\end{figure}
\begin{figure}[h] \centering
\includegraphics[width=7cm]{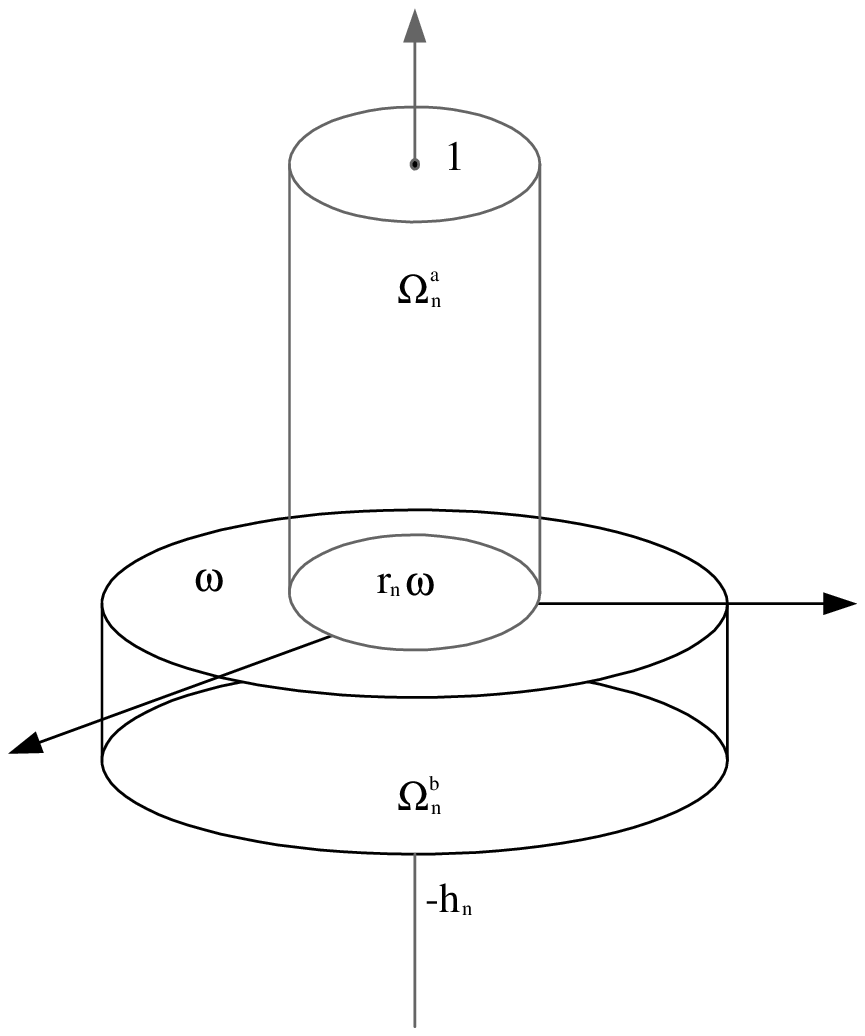}
\caption{the thin multidomain for $N=3$.} \label{domain1}
\end{figure}

In $\Omega_n$  consider the following eigenvalue problem:
\begin{equation}\label{opx}\left\{\begin{array}{ll}-\Delta U_n=\lambda U_n \hbox{ in
}\Omega_n,\\\\ U_n=0 \hbox{ on }
\Gamma_n=(r_n\omega\times\{1\})\cup(
\partial\omega\times ]-h_n,0[),\\\\\displaystyle{
\frac{\partial U_n}{\partial\nu}=0 \hbox{ on } \partial
\Omega_n\setminus\Gamma_n,}
\end{array}\right.
\end{equation}
where $\nu$  denotes the exterior
 unit normal to $\Omega_n$. Remark that $r_n\omega\times\{1\}$ is the top of the upper cylinder,
while $\partial\omega\times ]-h_n,0[$ is the lateral surface of
the second one.

It is well known (for instance, see Th. 6.2-1 in \cite{RT}) that,
for every $n \in \mathbb N$, there exists an increasing diverging
sequence of positive numbers $\{\lambda_{n,k}\}_{k \in \mathbb N}$
and a $L^2(\Omega_n)$-Hilbert orthonormal basis $\{U_{n,k}\}_{k
\in \mathbb N}$, such that $\{\lambda_{n,k}\}_{k \in \mathbb N}$
forms the set of all  the eigenvalues of Problem (\ref{opx}) and,
for every $k\in \mathbb N$,  $U_{n,k}\in {\cal V}_n=\left\{V\in
H^1(\Omega_n)\,:\,V=0 \hbox{ on }\Gamma_n\right\}$ is an
eigenvector of (\ref{opx}) with eigenvalue $\lambda_{n,k}$.
Moreover, $\left\{\lambda_{n,k}^{-\frac{1}{2}}U_{n,k}\right\}_{k
\in \mathbb N}$ is a ${\cal V}_n$-Hilbert orthonormal basis,  by
equipping ${\cal V}_n$
 with the inner product: $\displaystyle{(U,V)\in {\cal
V}_n\times{\cal V}_n\rightarrow}$ $\displaystyle{\int_{\Omega_n} D
UDVdx}$.

The aim of this paper is to study the asymptotic behavior of the
sequences $\{(\lambda_{n,k},U_{n,k}) \}_{n \in \mathbb N}$, for
every $k\in \mathbb N$,  as $n$ diverges.

We derive the limit eigenvalue problem (which is well posed in the
union of the limit domains: $]0,1[$ and $\omega$, with respective
dimension: $1$ and $N-1$) and the limit basis, and we precise how
these limits depend on the dimension $N$ and on limit
$\displaystyle{q=\lim_n\frac{h_n}{r_n^{N-1}}}$ of the ratio
between the volumes of the two cylinders.

In the sequel, $D$ stands for the gradient. Moreover,
$x=(x_1,\cdots , x_{N-1},x_N)=(x',x_N)$ denotes the generic point
of $\mathbb R^N$, and
 $D_{x'}$ and  $D_{x_N}$  stand for the gradient  with
respect to the first $N-1$ variables and for the derivative with
respect to the last variable, respectively. Furthermore, if $u$
depends only on one variable, $u'$ and $u''$ stand for the first
and second derivative, respectively.\smallskip

If $h_n \simeq r_n^{N-1}$, we obtain the following main result:

\begin{Theorem}\label{TH1}
Assume that
$$\displaystyle{\lim_n\frac{h_n}{r_n^{N-1}}=q\in]0,+\infty[}.$$

Then, there exists an increasing diverging sequence of positive
numbers $\{\lambda_{k}\}_{k \in \mathbb N}$ such that
\begin{equation}\label{ac}\lim_n\lambda_{{n},k}=\lambda_k, \quad\forall k \in \mathbb N,\end{equation}
and  $\{\lambda_{k}\}_{k \in \mathbb N}$ is the set of all the
eigenvalues of the following problem:
\begin{equation}\label{4bis}\left\{\begin{array}{l}-{u^a}''=\lambda
u^a\hbox{ in }]0,1[,\\\\
-\Delta {u^b}=\lambda u^b\hbox{ in }\omega,\\\\
u^a(1)=0,\quad{u^a}'(0)=0 ,\\\\
u^b=0\hbox{ on }\partial\omega,
\end{array}\right.
\end{equation}
if $N\geq3$; of the following one:
\begin{equation}\label{4}\left\{\begin{array}{l}-{u^a}''=\lambda
u^a\hbox{ in }]0,1[,\\\\
-{u^b}''=\lambda u^b\hbox{ in }]c,0[,\\\\
-{u^b}''=\lambda u^b\hbox{ in }]0,d[,\\\\
u^a(1)=0,\\\\
u^b(c)=0,\quad u^b(d)=0,\\\\
u^a(0)=u^b(0),\\\\
\vert\omega\vert {u^a}'(0)=q\left({u^b}'(0^-)-{u^b}'(0^+)\right),
\end{array}\right.
\end{equation}
if $N=2$ and $\omega=]c,d[$.

There exists an increasing sequence of positive integer numbers
$\{n_i\}_{i \in \mathbb N}$ and a sequence
$\{(u^a_{k},u^b_{k})\}_{k \in \mathbb N}\subset V$ (depending
possibly on the selected subsequence $\{n_i\}_{i \in \mathbb N}$),
where
$$V=\left\{(v^a,v^b)\in H^1(]0,1[)\times
H^1(\omega)\,:\,v^a(1)=0,\, v^b=0 \hbox{ on }
\partial\omega, (\hbox{and }v^a(0)=v^b(0)\hbox{ if
}N=2)\right\},$$ such that
\begin{equation}\label{5}\displaystyle{\lim_{i}\left(\int_{r_{n_i}\omega\times
]0,1[}\left\vert
U_{n_i,k}-\frac{u^a_k}{r_{n_i}^{\frac{N-1}{2}}}\right\vert^2+
\left\vert D_{x'}U_{n_i,k}\right\vert^2+\left\vert
\partial_{x_N}U_{n_i,k}-\frac{{u^a_k}'}{r_{n_i}^{\frac{N-1}{2}}}\right\vert^2dx\right)=0,}\end{equation}
\begin{equation}\label{6}\displaystyle{\lim_i\left(
\int_{\omega\times ]-h_{n_i},0[}\left\vert U_{n_i,k}-\frac{
q^{\frac{1}{2}}}{h_{n_i}^{\frac{1}{2}}}u^b_k\right\vert^2+
\left\vert D_{x'}U_{n_i,k}-\frac{
q^{\frac{1}{2}}}{h_{n_i}^{\frac{1}{2}}}Du^b_k\right\vert^2+\left\vert
\partial_{x_N} U_{n_i,k}\right\vert^2dx\right)=0.}\end{equation}
as $i\rightarrow+\infty$, for every $k\in  \mathbb N$, and
$u_{k}=(u^a_k,u^b_k)$ is an eigenvector of Problem (\ref{4bis}) if
$N\geq3$ (Problem (\ref{4}) if $N=2$) with eigenvalue
$\lambda_{k}$

 Moreover, $\{u_{k}\}_{k \in \mathbb N}$ is a $L^2(]0,1[)\times L^2(\omega)$-orthonormal basis
 with respect to the inner product: $ \displaystyle{\vert\omega\vert\int_0^1
u^a v^a dx_N+q\int_{\omega}u^b v^b dx'}$, and
$\{\lambda_k^{-\frac{1}{2}}u_{k}\}_{k\in \mathbb{N}}$ is a
$V$-Hilbert orthonormal basis with respect to the inner product:
$\displaystyle{\vert\omega\vert\int_0^1
{u^a}'{v^a}'dx_N+q\int_{\omega}D u^bDv^bdx'}$.
\end{Theorem}

If $N\geq3$, consider the following problems:
\begin{equation}\label{4'}\left\{\begin{array}{l}-{u^a}''=\lambda^a
u^a\hbox{ in }]0,1[,\\\\
u^a(1)=0,\quad
 {u^a}'(0)=0 ,
\end{array}\right.
\end{equation}
and
\begin{equation}\label{4"}\left\{\begin{array}{l}
-\Delta u^b=\lambda^b u^b\hbox{ in }\omega,\\\\
u^b=0\hbox{ on }\partial\omega.
\end{array}\right.
\end{equation}
By denoting with $\{\lambda^a_{k}\}_{k \in \mathbb N}$  the
increasing  sequence of all the eigenvalues of Problem (\ref{4'})
(i.e.
$\displaystyle{\left\{\left(\frac{\pi}{2}+k\pi\right)^2\right\}_{k\in
\mathbb N_0}}$, where $\mathbb N_0=\mathbb N\cup\{0\}$), and with
$\{\lambda^b_{k}\}_{k \in \mathbb N}$ the increasing
 sequence of all the eigenvalues of Problem (\ref{4"}), it is easily seen that:
$$\{\lambda_k :k\in \mathbb N\}=\{\lambda^a_k :k\in \mathbb N\}\cup\{\lambda^b_k :k\in \mathbb
N\}.$$ As regards as the multiplicity, if $\lambda$ is an
eigenvalue only of Problem (\ref{4'}), then it is a simple
eigenvalue of Problem (\ref{4bis}) with eigenvector $\left(\cos
\left(\frac{\pi}{2}+k\pi\right)x_N,0\right) $, for some
$k\in\mathbb N_0$. If $\lambda$ is an eigenvalue only of Problem
(\ref{4"}) with multiplicity $s$ (and with linearly independent
eigenvectors: $w^b_1,\cdots,w^b_s$), then it is an eigenvalue of
Problem (\ref{4bis}) with multiplicity $s$ and with linearly
independent eigenvectors: $ (0,w^b_1),\cdots,
 (0,w^b_s) $. If
$\lambda$ is an eigenvalue    of both Problem (\ref{4"}) with
multiplicity $s$ (and with linearly independent eigenvectors:
$w^b_1,\cdots,w^b_s$) and  Problem (\ref{4'}), then it is an
eigenvalue of Problem (\ref{4bis}) with multiplicity $s+1$ and
with linearly independent eigenvectors: $ (0,w^b_1),\cdots,
(0,w^b_s),\left(\cos \left(\frac{\pi}{2}+k\pi\right)x_N,0\right)$,
for some $k\in\mathbb N_0$. Roughly speaking, the eigenvalues of
Problem (\ref{4bis}) are obtained by gathering the eigenvalues of
Problem (\ref{4'}) and the eigenvalues of Problem (\ref{4"}).
Moreover, each eigenvalue preserves its multiplicity if it is
eigenvalue only of Problem (\ref{4'}) or only of Problem
(\ref{4"}); otherwise its multiplicity is obtained by adding the
multiplicity as eigenvalue
 of Problem (\ref{4'}) and the multiplicity as eigenvalue
of  Problem (\ref{4"}).

\noindent The limit $q$ of the ratio between the volumes of the
two cylinders does not intervene in the limit eigenvalue Problem
(\ref{4bis}). It appears  in the othonormal conditions:
\begin{equation}\label{"s"} \displaystyle{\vert\omega\vert\int_0^1 u^a_h u^a_k
dx_N+q\int_{\omega}u^b_h u^b_k dx'=\delta_{h,k}},\quad\forall
h,k\in \mathbb N,\end{equation}
\begin{equation}\label{"c"}\displaystyle{
\lambda_h^{-\frac{1}{2}}\lambda_k^{-\frac{1}{2}}\left(\vert\omega\vert\int_0^1{
u^a_h}'{u^a_k}'dx_N+q\int_{\omega}D u^b_hDu^b_kdx'
\right)=\delta_{h,k}},\quad\forall h,k\in \mathbb N,
\end{equation} where $\delta_{h,k}$ is the Kronecker's delta.
Moreover, it explicitly intervenes in  the corrector result
(\ref{6}), and, by  the previous othonormal conditions, also in
the corrector result (\ref{5}).

If $N=2$,  $q$ appears  also in the limit eigenvalue Problem
(\ref{4}). More precisely, if $N=2$,  the limit  problem in
$]0,1[$ is coupled with the limit problem in $\omega$ by the
junction conditions:
$$u^a(0)=u^b(0)\quad \hbox{ and } \quad\vert\omega\vert
{u^a}'(0)=q\left({u^b}'(0^-)-{u^b}'(0^+)\right).$$
 An easy computation shows that the
eigenvalues of Problem (\ref{4}) are given by  the positive
solutions of the following equation:
$$
\vert\omega\vert\sin(c\sqrt{\lambda})\sin(d\sqrt{\lambda})\cos(\sqrt{\lambda})-q\sin(\sqrt{\lambda})\sin(d\sqrt{\lambda})\cos(c\sqrt{\lambda})
+q\sin(\sqrt{\lambda})\sin(c\sqrt{\lambda})\cos(d\sqrt{\lambda})=0.$$
 For instance, if $\omega=]-1,1[$,  the set of the eigenvalues of Problem
(\ref{4}) is
$\displaystyle{\left\{\left(k\frac{\pi}{2}\right)^2\right\}_{k\in
\mathbb N}}$, and $\displaystyle{\left(k\frac{\pi}{2}\right)^2}$
is a simple eigenvalue if $k $ is odd, it is an eigenvalue with
multiplicity $2$ if $k$ is even. If $\omega=]-1,2[$, the set of
the eigenvalues of Problem (\ref{4}) is
$\displaystyle{\left\{\left(k{\pi}\right)^2\right\}_{k\in \mathbb
N}}\cup \displaystyle{\left\{\left(\pm\arccos\left(\pm
\sqrt{\frac{q}{4q+2\vert\omega\vert}}\right)+2k{\pi}\right)^2\right\}_{k\in
\mathbb N_0}}.$
\medskip

If  $h_n \ll r_n^{N-1}$, we obtain the following result:
\begin{Theorem}\label{TH2} Assume that
$$\displaystyle{\lim_n\frac{h_n}{r_n^{N-1}}=q=0}.$$

Then, there exists an increasing diverging sequence of positive
numbers $\{\lambda_{k}\}_{k \in \mathbb N}$ satisfying (\ref{ac}),
and $\{\lambda_{k}\}_{k \in \mathbb N}$ is the set of all the
eigenvalues of Problem (\ref{4bis}) if $N\geq3$, of the following
problem:
\begin{equation}\label{7}\left\{\begin{array}{l}-{u^a}''=\lambda
u^a\hbox{ in }]0,1[,\\\\
-{u^b}''=\lambda u^b\hbox{ in }]c,0[,\\\\
-{u^b}''=\lambda u^b\hbox{ in }]0,d[,\\\\
u^a(1)=0,\quad {u^a}'(0)=0,\\\\
u^b(c)=0,\quad u^b(0)=0, \quad u^b(d)=0,
\end{array}\right.\end{equation}
if $N=2$ and $\omega=]c,d[$.

 There exists an increasing sequence
of positive integer numbers $\{n_i\}_{i \in \mathbb N}$ and a
sequence $\{(u^a_{k},u^b_{k})\}_{k \in \mathbb N}\subset V_0$
(depending possibly on the selected subsequence $\{n_i\}_{i \in
\mathbb N}$), where
$$V_0=\left\{(v^a,v^b)\in H^1(]0,1[)\times
H^1(\omega)\,:\,v^a(1)=0,\, v^b=0 \hbox{ on }
\partial\omega, (\hbox{and }v^b(0)=0\hbox{ if }N=2)\right\},$$
satisfying (\ref{5}) and (\ref{6}) with $q=1$, and
$u_{k}=(u^a_k,u^b_k)$ is an eigenvector of Problem (\ref{4bis}) if
$N\geq3$ (Problem (\ref{7}) if $N=2$)
 with  eigenvalue $\lambda_{k}$.

 Moreover, $\{u_{k}\}_{k \in
\mathbb N}$ is a $L^2(]0,1[)\times L^2(\omega)$-orthonormal basis
 with respect to the inner product: $ \displaystyle{\vert\omega\vert\int_0^1
u^a v^a dx_N+\int_{\omega}u^b v^b dx'}$, and
$\{\lambda_k^{-\frac{1}{2}}u_{k}\}_{k\in \mathbb{N}}$ is a
$V_0$-Hilbert orthonormal basis with respect to the inner product:
$\displaystyle{\vert\omega\vert\int_0^1
{u^a}'{v^a}'dx_N+\int_{\omega}D u^bDv^bdx'}$.
\end{Theorem}

If $N\geq 3$,  one obtains the same result as in the previous case
when $q=1$ (see Theorem \ref{TH1}).

If $N=2$, consider the following  problems:
\begin{equation}\nonumber\left\{\begin{array}{l}-{u^a}''=\lambda^a
u^a\hbox{ in }]0,1[,\\\\
u^a(1)=0,\,\,\,
 {u^a}'(0)=0 ,
\end{array}\right.
\quad \left\{\begin{array}{l}-{u^b}''=\lambda^b
u^b\hbox{ in }]c,0[,\\\\
u^b(c)=0,\,\,\,
 {u^b}(0)=0 ,
\end{array}\right.
\quad \left\{\begin{array}{l}-{\overline
u^b}''=\overline{\lambda}^b \overline u^b\hbox{ in
}]0,d[,\\\\
\overline u^b(0)=0,\,\,\,
 {\overline
u^b}(d)=0 .
\end{array}\right.
\end{equation}
 Then, with  analogous considerations written
when $N\geq3$ after Theorem \ref{TH1}, it results that the
eigenvalues of Problem (\ref{7}) are obtained by gathering the
eigenvalues of these three problems, that is
$$\displaystyle{\{\lambda_k \,:\, k\in \mathbb N\}= \left\{\left(\frac{\pi}{2}+k\pi\right)^2\right\}_{k\in
\mathbb N_0}\cup\left\{\left(\frac{k\pi}{c}\right)^2\right\}_{k\in
\mathbb N}\cup\left\{\left(\frac{k\pi}{d}\right)^2\right\}_{k\in
\mathbb N},}$$ and that the multiplicity of each eigenvalue  is
less or equal than $3$. For instance, if $c=-1$ and $d=1$, the set
of the eigenvalues of Problem (\ref{7}) is
$\displaystyle{\left\{\left(k\frac{\pi}{2}\right)^2\right\}_{k\in
\mathbb N}}$, and $\displaystyle{\left(k\frac{\pi}{2}\right)^2}$
is a simple eigenvalue if $k $ is odd, it is an eigenvalue with
multiplicity $2$ if $k$ is even. If $c=-2$ and $d=2$, the set of
the eigenvalues of Problem (\ref{7}) is always
$\displaystyle{\left\{\left(k\frac{\pi}{2}\right)^2\right\}_{k\in
\mathbb N}}$, but  $\displaystyle{\left(k\frac{\pi}{2}\right)^2}$
is a simple eigenvalue if $k $ is even, it is an eigenvalue with
multiplicity $3$ if $k$ is odd.

\medskip

If $r_n^{N-1}\ll h_n$  we obtain the following partial result (the
result is complete for $N=2$) :

\begin{Theorem}\label{TH3} Assume that
\begin{equation}\label{ad}\displaystyle{\lim_n\frac{h_n}{r_n^{N-1}}=q=+\infty},\end{equation} and
\begin{equation}\label{ad1} h_n\ll -r_n^2\log r_n,\hbox{ if }N=3;\quad\quad h_n\ll
r_n^2, \hbox{ if } N\geq 4 .\end{equation}

Then, there exists an increasing diverging sequence of positive
numbers $\{\lambda_{k}\}_{k \in \mathbb N}$ satisfying (\ref{ac}),
and $\{\lambda_{k}\}_{k \in \mathbb N}$ is the set of all the
eigenvalues of Problem (\ref{4bis}) if $N\geq3$, of the following
problem:
\begin{equation}\label{70}\left\{\begin{array}{l}-{u^a}''=\lambda
u^a\hbox{ in }]0,1[,\\\\
-{u^b}''=\lambda u^b\hbox{ in }\omega=]c,d[,\\\\
u^a(1)=0,\quad u^a(0)=0,\\\\
u^b(c)=0, \quad u^b(d)=0,
\end{array}\right.\end{equation}
if $N=2$.

 There exists an increasing sequence of positive integer
numbers $\{n_i\}_{i \in \mathbb N}$ and a sequence
$\{(u^a_{k},u^b_{k})\}_{k \in \mathbb N}\subset V_\infty$
(depending possibly on the selected subsequence $\{n_i\}_{i \in
\mathbb N}$), where
$$V_\infty=\left\{(v^a,v^b)\in H^1(]0,1[)\times
H^1(\omega)\,:\,v^a(1)=0,\, v^b=0 \hbox{ on }
\partial\omega, (\hbox{and }v^a(0)=0\hbox{ if }N=2)\right\},$$
satisfying (\ref{5}) and (\ref{6}) with $q=1$,
 and $u_{k}=(u^a_k,u^b_k)$ is an eigenvector of Problem
(\ref{4bis}) if $N\geq3$  (Problem (\ref{70}) if $N=2$) with
eigenvalue $\lambda_{k}$.

 Moreover, $\{u_{k}\}_{k \in
\mathbb N}$ is a $L^2(]0,1[)\times L^2(\omega)$-orthonormal basis
 with respect to the inner product: $ \displaystyle{\vert\omega\vert\int_0^1
u^a v^a dx_N+\int_{\omega}u^b v^b dx'}$, and
$\{\lambda_k^{-\frac{1}{2}}u_{k}\}_{k\in \mathbb{N}}$ is a
$V_\infty$-Hilbert orthonormal basis with respect to the inner
product: $\displaystyle{\vert\omega\vert\int_0^1
{u^a}'{v^a}'dx_N+\int_{\omega}D u^bDv^bdx'}$.
\end{Theorem}

If $N\geq 3$, assumption (\ref{ad}) and (\ref{ad1}) reduce to
$$r_n^{2}\ll h_n\ll -r_n^2\log r_n,\hbox{ if }N=3;\quad\quad r_n^{N-1}\ll h_n\ll
r_n^2, \hbox{ if } N\geq 4 .$$

If $N\geq 3$,  one obtains the same result as in the previous
cases when $q=1$ or $q=0$ (see Theorem \ref{TH1} and Theorem
\ref{TH2}).

If $N=2$, by considering the  following problems:
\begin{equation}\nonumber\left\{\begin{array}{l}-{u^a}''=\lambda^a
u^a\hbox{ in }]0,1[,\\\\
u^a(1)=0,\quad
 {u^a}(0)=0 ,
\end{array}\right.
\quad \left\{\begin{array}{l}
-{u^b}''=\lambda^b u^b\hbox{ in }\omega=]c,d[,\\\\
u^b(c)=0,\quad u^b(d)=0,
\end{array}\right.
\end{equation}
 it is possible to repeat similarly the
analysis written when $N\geq3$ after Theorem \ref{TH1}. Then,
 the eigenvalues of Problem (\ref{70}) are obtained by gathering the
eigenvalues of these two last problems, that is
$$\displaystyle{\{\lambda_k \,:\, k\in \mathbb N\}=\left\{\left(k\pi\right)^2\right\}_{k\in \mathbb
N}\cup\left\{\left(\frac{k\pi}{d-c}\right)^2\right\}_{k\in \mathbb
N},}$$ and  the multiplicity of each eigenvalue is $1$ or $2$. For
instance, if $]c,d[=]-1,1[$, then  the set of all the eigenvalues
of Problem (\ref{70}) is
$\displaystyle{\left\{\left(k\frac{\pi}{2}\right)^2\right\}_{k\in
\mathbb N}}$, and  $\displaystyle{\left(k\frac{\pi}{2}\right)^2}$
is a simple eigenvalue if $k $ is odd, it is an eigenvalue with
multiplicity $2$ if $k$ is even. If
$\displaystyle{]c,d[=\left]-\frac{1}{2},\frac{1}{2}\right[}$, then
the set of all the eigenvalues of Problem (\ref{70}) is
$\left\{\left(k\pi\right)^2\right\}_{k\in \mathbb N}$ and each
eigenvalue has multiplicity $2$. If
$\displaystyle{]c,d[=\left]-\frac{\pi}{2},\frac{\pi}{2}\right[}$,
then the set of all the eigenvalues of Problem (\ref{70}) is
$\{1,4,9,\pi^2,16,25,36, 4\pi^2, 49, 64, 81,9\pi^2,\cdots\}$ and
each eigenvalue has multiplicity $1$.\medskip

In summarizing:

$\bullet$ the limit eigenvalue problem reduces to a problem in
$]0,1[$ and a problem in $\omega$ (recall that $\omega$ has
dimension dimension $N-1$).

$\bullet$ if $N=2$, the limit eigenvalue problem
  depends on $q$ and it is  coupled if $q\in ]0,+\infty[$, uncoupled if $q=0$ or $q=+\infty$. Also the limit
eigenvector basis (and the orthonormal limit conditions) depends
on $q$.

$\bullet$ If $N\geq 3$ the limit eigenvalue problem  is
independent of $q$ and it is uncoupled (at least in the considered
cases, i.e. for $h_n\ll -r_n^2\log r_n$, if $N=3$; for $ h_n\ll
r_n^2$, if $N\geq 4 $).  The limit eigenvector basis (and the
orthonormal limit conditions)  depend on $q$. The orthonormal
limit conditions coincide for $q\in\{0,1,+\infty\}$.

$\bullet$ If the limit eigenvalue problem is uncoupled (that is
when $N=2$ and $q\in\{0,+\infty\}$, or when $N\geq 3$), the limit
 eigenvalues are obtained by gathering the eigenvalues
of a Laplace problem in $]0,1[$ and the eigenvalues of one or two
Laplace problems in a domain of dimension $N-1$ and by "adding the
multiplicities".

Remark that, in all the previous theorems, by virtue of
(\ref{ac}), the multiplicity of $\lambda_{n,k}$, for $n$ large
enough, is less or equal  than the multiplicity of $\lambda_{k}$.
Consequently, if $\lambda_{k}$ is simple, also $\lambda_{n,k}$ is
simple, for $n$ large enough. Then, by arguing as in \cite{V} (see
also \cite{CSJP}), if $\lambda_{k}$ is simple, by fixing one of
the two normalized eigenvectors $u_k$ of the limit problem with
eigenvalue $\lambda_{k}$, it is possible to choose, for $n$ large
enough, one of the two normalized eigenvectors $U_{n,k}\in {\cal
V}_n$  of Problem (\ref{opx}) with eigenvalue $\lambda_{n,k}$ such
that convergences (\ref{5}) and (\ref{6}) hold true for the whole
sequence.

Point out that it is not necessary that the two cylinders are
scaled to the same one or that the first cylinder has height $1$.
In fact, the results do not essentially change if one assumes
$\Omega_n=\left(r_n\omega_a\times [0,l[\right)\bigcup$
$\left(\omega_b\times ]-h_n,0[\right)$, with
$0'\in\omega^a\subset\omega^b$ and $l\in]0,+\infty[$.
\medskip

 In Section \ref{mainresult}, after having
reformulated the problem on a fixed domain through appropriate
rescalings of the kind proposed by P.G. Ciarlet and P. Destuynder
in \cite{CD},  and having introduced suitable weighted inner
products, by using the min-max Principle we obtain {\em a priori}
estimates (with respect to $n$) for the sequences
$\{\lambda_{n,k}\}_{k \in \mathbb N}$ (see Proposition
\ref{stimeautovalori}). Then, by making use of the method of
oscillating
    test functions,
    introduced by L. Tartar in \cite{T}, by applying some results obtained  by A. Gaudiello, B. Gustafsson, C. Lefter and J.
    Mossino in \cite{GGLM1} and \cite{GGLM2}
    and  by adapting the techniques used by M. Vanninathan in \cite{V},
    we derive the limit eigenvalue problem and the limit of the rescaled basis, as $n\rightarrow+\infty$,  in the case $h_n \simeq r_n^{N-1}$
    (see Theorem \ref{maintheorem} and the proof of Theorem \ref{TH1} at the end of Section
    \ref{mainresult}). The cases  $h_n \ll r_n^{N-1}$ and  $ r_n^{N-1}\ll h_n
    $ are sketched in Section \ref{Sec3} and Section \ref{Sec4},
    respectively.\medskip

For  the study of thin  multi-structures we refer to \cite{C1},
\cite{CD}, \cite{CSJP}, \cite{G}, \cite{KMM}, \cite{LeD}, \cite{P}
and the references quoted therein. For a thin multi-structure as
considered in this paper, we refer to \cite{GGLM0}, \cite{GGLM1},
\cite{GGLM2},  \cite{GMMMS1}, \cite{GMMMS2} and \cite{GZ}. For the
study of the spectrum of the Laplace operator in a thin tube with
a Dirichlet condition on its boundary we refer to  \cite{BMT}. For
the study of the homogenization of the spectrum of the Laplace
operator in a periodic perforated domain with different boundary
conditions on the holes we refer to \cite{V}.

\section{The case $\displaystyle{\lim_n\frac{h_n}{r_n^{N-1}}=q\in]0,+\infty[}$}\label{mainresult}

In the sequel, $\Omega^a = \omega\times ]0,1[$, $\Omega^b =
\omega\times ]-1,0[$, and $\left\{r_n\right\}_{n \in \mathbb N}$
and $\left\{h_n\right\}_{n \in \mathbb N}\subset]0,1[$ are two
sequences such that
\begin{equation}\label{hrzero}
\lim_{n }h_n=0= \lim_{n }r_n.
\end{equation}

For every $n \in \mathbb N$, let $H_n$ be the space
$L^2(\Omega^a)\times L^2(\Omega^b) $ equipped
 with the inner product:
\begin{equation}\label{prodottoscalaren} \begin{array}{l}
(\cdot,\cdot)_n:(u,v)=((u^a,u^b),
(v^a,v^b))\in\left(L^2(\Omega^a)\times L^2(\Omega^b)\right) \times
\left(L^2(\Omega^a)\times L^2(\Omega^b)\right)
\longrightarrow\\\\\displaystyle{(u,v)_n=\int_{\Omega^a} u^a v^a
dx+\frac{h_n}{r_n^{N-1}}\int_{\Omega^b}u^b v^b
dx,}\end{array}\end{equation} and let  $V_n$ be the  space:
\begin{equation}\label{Vn}\begin{array}{l}
\displaystyle{ \big\{v=(v^a, v^b) \in H^1 (\Omega^a)\times H^1
(\Omega^b)\,: \, v^a=0 \hbox{ on
}\omega\times \{1\},   }\\\\
 \quad\quad\quad\quad\quad\quad\quad\displaystyle {v^b=0 \hbox{ on } \partial\omega\times
]-1,0[,\,\, v^a(x',0)= v^b(r_nx',0)\hbox{ for }x' \hbox{ a.e. in
}\omega \big\}}\end{array}\end{equation} equipped
 with the inner product:
\begin{equation}\label{formabilinearen}\begin{array}{l}
a_n:(u,v)=((u^a,u^b), (v^a,v^b))\in V_n\times V_n\longrightarrow a_n(u,v)=\\\\
\displaystyle{\int_{\Omega^a}\frac{1}{r_n^2}D_{x'}
u^aD_{x'}v^a+\partial_{x_N}
u^a\partial_{x_N}v^adx+\frac{h_n}{r_n^{N-1}}\int_{\Omega^b}D_{x'}
u^bD_{x'}v^b+\frac{1}{h_n^2}\partial_{x_N}
u^b\partial_{x_N}v^bdx.}
\end{array}\end{equation}
The choice of $H_n$ and  $V_n$  will be justified in the proof of
Theorem \ref{TH1}, at the end of this section. Now, remark that
the norm induced on $V_n$ by the inner product
 $a_n(\cdot,\cdot)$ is equivalent to the $\left(H^1(\Omega^a)\times H^1(\Omega^b)\right)$-norm, and
  the norm induced on $H_n$ by the inner product
 $(\cdot,\cdot)_n$ is equivalent to the $\left(L^2(\Omega^a)\times
 L^2(\Omega^b)\right)$-norm. Consequently, $V_n$ is continuously and
compactly embedded into $H_n$. Moreover, since
$C_0^\infty(\Omega^a)\times C_0^\infty(\Omega^b)\subset V_n$, it
results that $V_n$ is dense in $H_n$. Then, (for instance, see Th.
6.2-1 and Th. 6.2-2 in \cite{RT})  for every $n \in \mathbb N$,
there exists an increasing diverging sequence of positive numbers
$\{\lambda_{n,k}\}_{k \in \mathbb N}$ and a $H_n$-Hilbert
orthonormal basis $\{u_{n,k}\}_{k \in \mathbb N}$, such that
$\{\lambda_{n,k}\}_{k \in \mathbb N}$ forms the set of all  the
eigenvalues of the following problem:
\begin{equation}\label{eigenproblem}\left\{\begin{array}{l}
\displaystyle{u_n\in { V}_n,}\\\\
\displaystyle{a_n(u_n, v)=\lambda(u_n, v)_n,\quad \forall v\in {
V}_n,}
\end{array}\right.\end{equation}
and, for every $k\in \mathbb N$,  $u_{n,k}\in V_n$ is an
eigenvector of (\ref{eigenproblem}) with eigenvalue
$\lambda_{n,k}$. Moreover,
$\left\{\lambda_{n,k}^{-\frac{1}{2}}u_{n,k}\right\}_{k \in \mathbb
N}$ is a $V_n$-Hilbert orthonormal basis. Furthermore, for every
$k\in \mathbb N$, $\lambda_{n,k}$ is characterized by the
following min-max Principle:
\begin{equation}\displaystyle{\label{minmax}\lambda_{n,k}=\min_{{\cal E}_k\in{\cal F}_k}\max_{v\in{\cal E}_k,
\,\,v\neq0}\frac{a_n(v,v)}{(v,v)_n},}\end{equation}  where ${\cal
F}_k$ is the set of the subspaces ${\cal E}_k$ of $V_n$ with
dimension $k$.

The min-max Principle provides the following {\em a priori}
estimates for the eigenvalues of  Problem (\ref{eigenproblem}):
\begin{Proposition}\label{stimeautovalori}
For every $n ,k\in \mathbb N$, let $\lambda_{n,k}$ be as above.
Then, it results that
\begin{equation}\label{stima}0<\lambda_{n,k}\leq 2^kk^2\pi^2,\quad\forall n, k \in \mathbb N.\end{equation}
\end{Proposition}
\begin{proof}[Proof]
It is well known that $\left\{j^2\pi^2\right\}_{j \in \mathbb N}$
is the sequence of the eigenvalues of the following problem:
\begin{equation}\label{problemaausiliare}\left\{\begin{array}{l}
\displaystyle{-y''(t)=\lambda y(t)\hbox{ in }]0,1[,}\\\\
y(1)=0= y(0),
\end{array}\right.\end{equation}
and $\{\sqrt{2}\sin(j\pi t)\}_{j \in \mathbb N}$ is a
$L^2(]0,1[)$-Hilbert orthonormal basis such that,  for every $j\in
\mathbb N$, $\sqrt{2}\sin(j\pi t)$ is an eigenvector of
(\ref{problemaausiliare}) with eigenvalue $j^2\pi^2$.

For every $j \in \mathbb N$, set
\begin{equation}\nonumber \zeta_j(x',x_N)=\left\{\begin{array}{l}y_{j}(x_N)=\sqrt{2}\sin(j\pi x_N),\hbox{
if
}(x',x_N)\in \Omega^a,\\\\
0,\hbox{ if }(x',x_N)\in\Omega^b.
\end{array}\right.\end{equation}

Fix $k \in \mathbb N$ and  set
$\displaystyle{Z_k=\left\{\sum_{j=1}^k\alpha_j \zeta_j :
\alpha_1,\cdots,\alpha_k\in \mathbb R\right\}}$. Then, for every
$n \in \mathbb N$,  $Z_k$ is a subspace of $V_n$  of dimension
$k$. Consequently, by applying the min-max Principle
(\ref{minmax}), it results that:
\begin{equation}\nonumber\begin{array}{l}
\displaystyle{\lambda_{n,k}\leq\max_{\zeta\in
Z_k-\{0\}}\frac{a_n(\zeta,\zeta)}{(\zeta,\zeta)_n}=
\max_{(\alpha_1,\cdots,\alpha_k)\in  \mathbb
R^k-\{0\}}\frac{\displaystyle{\int_0^1\left\vert\sum_{j=1}^k\alpha_jy'_j(t)
\right\vert^2dt}}
{\displaystyle{\int_0^1\left\vert\sum_{j=1}^k\alpha_jy_{j}(t)\right\vert^2dt}}\leq}\\\\
\displaystyle{\max_{(\alpha_1,\cdots,\alpha_k)\in  \mathbb
R^k-\{0\}}\frac{2^k\displaystyle{\sum_{j=1}^k\alpha_j^2\int_0^1\left\vert
y'_j(t)\right\vert^2dt}} {\displaystyle{\sum_{j=1}^k\alpha_j^2}}=
\max_{(\alpha_1,\cdots,\alpha_k)\in  \mathbb
R^k-\{0\}}\frac{2^k\displaystyle{\sum_{j=1}^k\alpha_j^2j^2\pi^2}}{\displaystyle{
\sum_{j=1}^k\alpha_j^2}}\leq}\\\\
\displaystyle{\max_{(\alpha_1,\cdots,\alpha_k)\in  \mathbb
R^k-\{0\}}\frac{2^k
k^2\pi^2\displaystyle{\sum_{j=1}^k\alpha_j^2}}{\displaystyle{\sum_{j=1}^k\alpha_j^2}}=2^k
k^2\pi^2,\quad\forall n\in\mathbb N.}
\end{array}\end{equation}
\end{proof}
\begin{Remark}\label{remarkstimaautov} By examining the previous proof, it is evident that
estimate (\ref{stima}) holds again true, if one replaces
$(\cdot,\cdot)_n$ and $a_n(\cdot,\cdot)$ with $c_n(\cdot,\cdot)_n$
and $c_na_n(\cdot,\cdot)$, respectively, where $\{c_n\}_{n \in
\mathbb N}\subset ]0,+\infty[$. \end{Remark}

By using a diagonal argument,  the following result is an
immediate consequence of Proposition \ref{stimeautovalori}.
\begin{Corollary}\label{corollariostimeautovalori}
For every $n ,k\in \mathbb N$, let $\lambda_{n,k}$ be as above.
Then, there exists an increasing sequence of positive integer
numbers $\{n_i\}_{i \in \mathbb N}$ and an increasing sequence of
no negative  numbers $\{\lambda_{k}\}_{k \in \mathbb N}$, such
that
$$\lim_i\lambda_{{n_i},k}=\lambda_k, \quad\forall k \in \mathbb N.$$
\end{Corollary}
\begin{Remark}In the sequel it will be shown that $\{\lambda_{k}\}_{k
\in \mathbb N}$ is a sequence of positive numbers and the
convergence holds true for the whole sequence
$\{\lambda_{n,k}\}_{n \in \mathbb N}$.
\end{Remark}

To characterize the sequence $\{\lambda_{k}\}_{k \in \mathbb N}$,
 consider the following limit:
\begin{equation}\label{volumiuguali}\lim_n\frac{h_n}{r_n^{N-1}}=q,\end{equation}
which exists always for  a subsequence. If
\begin{equation}\label{qfinitomaggioredizero}q\in ]0,+\infty[,\end{equation}
in the  space:
\begin{equation}\label{H}H=\{v=(v^a, v^b) \in L^2 (\Omega^a)\times
L^2 (\Omega^b): v^a \hbox{ is indepenendent of }x',\, v^b \hbox{
is indep. of }x_N \}
\end{equation}
 introduce the inner product:
\begin{equation}\label{prodottoscalarelim} \begin{array}{l}
[\cdot,\cdot]_{q}:(u,v)=((u^a,u^b), (v^a,v^b))\in H \times
H\longrightarrow \displaystyle{\vert\omega\vert\int_0^1 u^a v^a
dx_N+q\int_{\omega}u^b v^b dx',}\end{array}\end{equation} and in
the space:
\begin{equation}\label{limitspace}V=\left\{\begin{array}{ll}\begin{array}{l}
\Big\{v=(v^a, v^b) \in H^1 (\Omega^a)\times H^1 (\Omega^b):\\\\
v^a \hbox{ is independent of }x',\, v^b \hbox{ is independent of
}x_N,\\\\ v^a(1)=0,\,\, v^b=0 \hbox{ on } \partial\omega,\\\\
v^a(0)=v^b(0')\Big\},\end{array}&\quad\quad\hbox{ if }  N=2, \\\\\\
\begin{array}{l}\Big\{v=(v^a, v^b) \in H^1 (\Omega^a)\times H^1
(\Omega^b):\\\\ v^a \hbox{ is
independent of }x',\, v^b \hbox{ is independent of }x_N,\\\\
v^a(1)=0,\,\, v^b=0 \hbox{ on } \partial\omega \Big\},\end{array}
&\quad\quad\hbox{ if } 3\leq N
\end{array}\right.\end{equation}  introduce the inner product:
\begin{equation}\label{formabilineare}\begin{array}{l}
{\alpha}_{q}:(u,v)=((u^a,u^b), (v^a,v^b))\in V\times
V\longrightarrow {\alpha}_{q}(u,v)=\\\\
\displaystyle{\vert\omega\vert\int_0^1\partial_{x_N}
u^a\partial_{x_N}v^adx_N+q\int_{\omega}D_{x'} u^bD_{x'}v^bdx' .}
\end{array}\end{equation}
Remark that the norm induced on $V$ by the inner product
 ${\alpha}_{q}(\cdot,\cdot)$ is equivalent to the $\left(H^1(]0,1[)\times H^1(\omega)\right)$-norm, and
  the norm induced on $H$ by the inner product
 $[\cdot,\cdot]_q$ is equivalent to the $\left(L^2(]0,1[)\times
 L^2(\omega)\right)$-norm. Consequently, $V$ is continuously and
compactly embedded into $H$. Moreover, since
$C_0^\infty(]0,1[)\times C_0^\infty(\omega)\subset V$ if $N\geq 3$
($C_0^\infty(]0,1[)\times\{v\in C_0^\infty(\omega):v(0)=0\}\subset
V$ if $N=2$), it results that $V$ is dense in $H$. Then, one can
consider the following eigenvalue problem:
\begin{equation}\label{limiteigenproblem}\left\{\begin{array}{l}
\displaystyle{u\in  V,}\\\\
\displaystyle{\alpha_q(u, v)=\lambda[u, v]_{q},\quad \forall v\in
{ V}},
\end{array}\right.\end{equation}
for which all the classic results hold true (see Th. 6.2-1 and
6.2-2 in \cite{RT}).

The main result of this paper is  the following one:
\begin{Theorem}\label{maintheorem}
For every $n \in \mathbb N$, let $\{\lambda_{n,k}\}_{k \in \mathbb
N}$ be an increasing diverging sequence of  all the eigenvalues of
Problem (\ref{prodottoscalaren})$\div$(\ref{eigenproblem}), and
$\{u_{n,k}\}_{k \in \mathbb N}$ be a $H_n$-Hilbert orthonormal
basis such that
$\left\{\lambda_{n,k}^{-\frac{1}{2}}u_{n,k}\right\}_{k \in \mathbb
N}$ is a $V_n$-Hilbert orthonormal basis and, for every $k\in
\mathbb N$, $u_{n,k}$ is an eigenvector of Problem
(\ref{prodottoscalaren})$\div$(\ref{eigenproblem}) with eigenvalue
$\lambda_{n,k}$.   Assume (\ref{hrzero}), (\ref{volumiuguali}) and
(\ref{qfinitomaggioredizero}).

Then, there exists an increasing diverging sequence of positive
numbers $\{\lambda_{k}\}_{k \in \mathbb N}$ such that
$$\lim_n\lambda_{{n},k}=\lambda_k, \quad\forall k \in \mathbb N,$$
and $\{\lambda_{k}\}_{k \in \mathbb N}$ is the set of all  the
eigenvalues of Problem (\ref{H})$\div$(\ref{limiteigenproblem}).
 Moreover,
there exists an increasing sequence of positive integer numbers
$\{n_i\}_{i \in \mathbb N}$ and a $(H,[\cdot,\cdot]_q)$-Hilbert
orthonormal basis $\{u_{k}\}_{k \in \mathbb N}$ (depending
possibly on the selected subsequence $\{n_i\}_{i \in \mathbb N}$)
such that, for every $k\in \mathbb N$, $u_{k}\in V$ is an
eigenvector of Problem (\ref{H})$\div$(\ref{limiteigenproblem})
with eigenvalue $\lambda_{k}$, and
\begin{equation}\label{1} u_{n_i,k}\rightarrow u_k\hbox{ strongly in
}H^1(\Omega^a)\times H^1(\Omega^b), \quad\forall k \in \mathbb N,
\end{equation}
as $i\rightarrow+\infty$,
\begin{equation}\label{2} \frac{1}{r_{n}}D_{x'}u_{n,k}^a\rightarrow 0'\hbox{ strongly in
}(L^2(\Omega^a))^{N-1},\quad\forall k \in \mathbb N,
\end{equation}
\begin{equation}\label{3}
\frac{1}{h_{n}}\partial_{x_N}u_{n,k}^b\rightarrow 0\hbox{ strongly
in }L^2(\Omega^b), \quad\forall k \in \mathbb N,
\end{equation}
as $n\rightarrow+\infty$. Furthermore,
$\{\lambda_k^{-\frac{1}{2}}u_{k}\}_{k\in \mathbb{N}}$ is a
$(V,\alpha_q)$-Hilbert orthonormal basis.
\end{Theorem}
\begin{proof}[Proof]$\,$ The proof will be performed in three
steps.

\noindent $\underline{\hbox {Step 1.}}$ This step is devoted to
derive the limit eigenvalue problem.

 It results that
\begin{equation}\label{eigenproblemperun,k}\left\{\begin{array}{l}
\displaystyle{u_{n,k}\in { V}_n,}\\\\
\displaystyle{a_n(u_{n,k}, v)=\lambda_{n,k}(u_{n,k}, v)_n,\quad
\forall v\in { V}_n,}
\end{array}\right.\quad\forall n,k \in \mathbb
N,\end{equation}
\begin{equation}\label{deltahk}
(u_{n,k}, u_{n,h})_n=\delta_{h,k},\quad \forall n,k,h, \in \mathbb
N,
\end{equation}
where
\begin{equation}\label{delta}\delta_{h,k}=\left\{\begin{array}{l}1,\hbox{
if }h=k,\\0,\hbox{ if }h\not= k.\end{array}\right.\end{equation}

By choosing $v=u_{n,k}$ in (\ref{eigenproblemperun,k}), and by
taking into account (\ref{deltahk}) and Proposition
\ref{stimeautovalori}, it follows that
\begin{equation}\nonumber
a_n(u_{n,k}, u_{n,k})=\lambda_{n,k}(u_{n,k},
u_{n,k})_n=\lambda_{n,k}\leq 2^kk^2\pi^2,\quad\forall n ,k \in
\mathbb N.
\end{equation}
Consequently, by virtue of  definition (\ref{formabilinearen}) of
$a_n$ and of  assumptions (\ref{hrzero}), (\ref{volumiuguali}) and
(\ref{qfinitomaggioredizero}), by applying Proposition 2.1 in
\cite{GGLM1} (which permits to obtain the junction condition when
$N=2$) and by using a diagonal argument, it is easily seen that
there exists an increasing sequence of positive integer numbers
$\{n_i\}_{i \in \mathbb N}$, an increasing sequence of no negative
numbers $\{\lambda_{k}\}_{k \in \mathbb N}$ and a sequence
$\left\{u_{k}\right\}_{k \in \mathbb N}\subset V$ (depending
possibly on the selected subsequence $\{n_i\}_{i \in \mathbb N}$)
such that
\begin{equation}\label{convlambda}\lim_i\lambda_{{n_i},k}=\lambda_k, \quad\forall k \in \mathbb N,\end{equation}
\begin{equation}\label{limitedeboleautovettori}u_{n_i,k}\rightharpoonup u_k\hbox{ weakly in
}H^1(\Omega^a)\times H^1(\Omega^b)\hbox{ and strongly in
}L^2(\Omega^a)\times L^2(\Omega^b), \quad\forall k \in \mathbb N,
\end{equation}
as $i\rightarrow +\infty$.

By passing to the limit in (\ref{eigenproblemperun,k}) and
(\ref{deltahk}), as $n_i\rightarrow +\infty$, by making use of
(\ref{volumiuguali}), (\ref{qfinitomaggioredizero}),
(\ref{convlambda}), (\ref{limitedeboleautovettori}),  and by
applying the Dirichlet version of Theorem 1.1 in  \cite{GGLM2}
(see also Remark 1.4 in \cite{GGLM2}), with $f=\lambda_k u_k$, it
turns out that (\ref{1})$\div$(\ref{3}) hold true and
\begin{equation}\label{problemalimite}\left\{\begin{array}{l}
\displaystyle{u_{k}\in { V},}\\\\
\displaystyle{\alpha_q(u_{k}, v)=\lambda_{k}[u_{k}, v]_q,\quad
\forall v\in { V},}
\end{array}\right.\quad\forall k\in \mathbb N,\end{equation}
\begin{equation}\label{deltahklim}
[u_{k}, u_{h}]_q=\delta_{h,k},\quad\forall k,h\in \mathbb N.
\end{equation}

In particular, (\ref{problemalimite}) and (\ref{deltahklim})
provide that
\begin{equation}\nonumber
\alpha_q(\lambda_{k}^{-\frac{1}{2}}u_{k},\lambda_{h}^{-\frac{1}{2}}
u_{h})=\delta_{h,k},\quad\forall k,h\in \mathbb N,
\end{equation}
Hence,  $\{\lambda_{k}\}_{k \in \mathbb N}$ is a sequence of
 positive numbers. Moreover, remark  that
\begin{equation}\label{lomautovalori}
\displaystyle{\lim_{k}\lambda_{k}=+\infty.}
\end{equation}
In fact, or (\ref{lomautovalori}) holds true, or
$\{\lambda_{k}\}_{k \in \mathbb N}$ is a finite set. In the second
case, by virtue of (\ref{deltahklim}),  Problem
(\ref{problemalimite}) would admit an eigenvalue of infinite
multiplicity. But this is not possible, due to  the Fredholm's
alternative Theorem.

\noindent $\underline{\hbox {Step 2.}}$
 This step is devoted to prove that there not exist $(\overline
u, \overline \lambda)\in V\times\mathbb{R}$ satisfying the
following problem:
\begin{equation}\label{probleoverlinelambda}\left\{\begin{array}{l}
\displaystyle{\overline
u\in { V},}\\\\
\displaystyle{\alpha_q(\overline u, v)=\overline\lambda [\overline
u, v]_q,\quad
\forall v\in { V},}\\\\

[\overline u, u_k]_q=0,\quad\forall k\in\mathbb N\\\\

[\overline u,\overline u ]_q=1.
\end{array}\right.\end{equation}
To this aim, the proof of  Theorem  9.2 in  \cite{V} (see also
\cite{CSJP}) will be
 adapted to our case.

By arguing by contradiction, assume that there exists $(\overline
u, \overline \lambda)\in V\times\mathbb{R}$ satisfying Problem
(\ref{probleoverlinelambda}).

First, remark that (\ref{lomautovalori}) provides the existence of
 ${\overline k}\in \mathbb{N}$ such that
\begin{equation}\label{maggiore}
\overline\lambda< \lambda_{\overline k}.
\end{equation}

For every $n\in \mathbb{N}$, let $\varphi_n$ be the solution of
the following problem:
\begin{equation}\label{problema varphin}\left\{\begin{array}{l}
\displaystyle{\varphi_n\in { V}_n,}\\\\
\displaystyle{a_n(\varphi_n, v)=\overline\lambda(\overline u,
v)_n,\quad \forall v\in { V}_n.}
\end{array}\right.\end{equation}
Then,   the Dirichlet version of Theorem 1.1 in \cite{GGLM2} (see
also Remark 1.4 in \cite{GGLM2}), with $f=\overline\lambda
\overline u $,  entails the existence of $\varphi\in V$ such that
\begin{equation}\label{limitedebolevarphin}\varphi_n\rightharpoonup \varphi\hbox{ weakly in
}H^1(\Omega^a)\times H^1(\Omega^b)\hbox{ and strongly in
}L^2(\Omega^a)\times L^2(\Omega^b),
\end{equation}
as $n\rightarrow +\infty$, and
\begin{equation}\label{ccc}\left\{\begin{array}{l}
\displaystyle{\varphi\in V,}\\\\
\displaystyle{\alpha_q(\varphi, v)=\overline\lambda[\overline u,
v]_q,\quad \forall v\in V.}
\end{array}\right.\end{equation}
By comparing (\ref{probleoverlinelambda}) with (\ref{ccc}) and by
virtue of the uniqueness of $\varphi$, it turns out that
\begin{equation}\label{uguale}\varphi=\overline u.
\end{equation}

For every $n\in \mathbb{N}$, set
$$v_n=\varphi_n-\sum_{i=1}^{\overline k}
(\varphi_n, u_{n,i})_n u_{n,i}\in V_n.$$ Then,  by virtue of
(\ref{deltahk}), it results that
\begin{equation}\nonumber\begin{array}{l}\displaystyle{
(v_n, u_{n,k})_n= (\varphi_n, u_{n,k})_n-\sum_{i=1}^{\overline k}
(\varphi_n, u_{n,i})_n (u_{n,i}, u_{n,k})_n =(\varphi_n,
u_{n,k})_n- (\varphi_n, u_{n,k})_n
=0,}\\\\\displaystyle{\quad\quad\quad\quad\quad\quad\quad\quad\quad\quad\quad\quad
\quad\quad\quad\quad\quad\quad\quad\quad\quad\quad\quad\quad\forall
k\in\{1,\cdots,\overline k\},\quad\forall n\in \mathbb{N}.}
\end{array}\end{equation}
Consequently, by recalling that
$\left\{\lambda_{n,k}^{-\frac{1}{2}}u_{n,k}\right\}_{k \in \mathbb
N}$ is a $V_n$-Hilbert orthonormal basis, that $u_{n,k}$ solves
Problem (\ref{eigenproblemperun,k}), that $\{\lambda_{n,k}\}_{k
\in \mathbb N}$ is an increasing sequence and that $\{u_{n,k}\}_{k
\in \mathbb N}$ is a $H_n$-Hilbert orthonormal basis,
   it follows that
\begin{equation}\label{classicalpeoperties}\begin{array}{l}\displaystyle{
a_n(v_n,v_n)= \sum_{k=1}^{+\infty}\vert
a_n(\lambda_{n,k}^{-\frac{1}{2}}u_{n,k},v_n)\vert^2=\sum_{k=1}^{+\infty}\lambda_{n,
k}^{-1}\lambda_{n, k}^2(u_{n,k},v_n)_n^2=}\\\\\displaystyle{
 \sum_{k=\overline k +1}^{+\infty}\lambda_{n,
k}(u_{n,k},v_n)_n^2 \geq\lambda_{n,\overline k+1}\sum_{k=\overline
k +1}^{+\infty}(u_{n,k},v_n)_n^2 =\lambda_{n,\overline
k+1}(v_n,v_n)_n,\quad\forall n\in \mathbb{N}.}
\end{array}\end{equation}

In what concerns the first term in (\ref{classicalpeoperties}),
\begin{equation}\nonumber\begin{array}{l}\displaystyle{
a_n(v_n,v_n)=a_n\left(\varphi_n-\sum_{k=1}^{\overline k}
(\varphi_n, u_{n,k})_n u_{n,k},\varphi_n-\sum_{j=1}^{\overline k}
(\varphi_n,
u_{n,j})_n u_{n,j}\right)= }\\\\
\displaystyle{a_n(\varphi_n,\varphi_n)-2\sum_{k=1}^{\overline
k}(\varphi_n, u_{n,k})_n
a_n(\varphi_n,u_{n,k})+\sum_{k,j=1}^{\overline k}(\varphi_n,
u_{n,k})_n (\varphi_n, u_{n,j})_n
a_n(u_{n,k},u_{n,j}),\quad\forall n\in \mathbb{N}.}
\end{array}\end{equation}
from which, by virtue  of (\ref{problema varphin}),
(\ref{eigenproblemperun,k}) and (\ref{deltahk}),  it follows that
\begin{equation}\nonumber\begin{array}{l}\displaystyle{
a_n(v_n,v_n)= \overline\lambda(\overline u,\varphi_n)_n
-2\overline\lambda\sum_{k=1}^{\overline k}(\varphi_n, u_{n,k})_n
(\overline u,u_{n,k})_n +\sum_{k,j=1}^{\overline k}(\varphi_n,
u_{n,k})_n (\varphi_n, u_{n,j})_n
\lambda_{n,k}(u_{n,k},u_{n,j})_n=} \\\\\displaystyle{
\overline\lambda(\overline u,\varphi_n)_n
-2\overline\lambda\sum_{k=1}^{\overline k}(\varphi_n, u_{n,k})_n
(\overline u,u_{n,k})_n +\sum_{k=1}^{\overline k}(\varphi_n,
u_{n,k})_n^2 \lambda_{n,k},\quad\forall n\in \mathbb{N}.}
\end{array}\end{equation}
Consequently, by virtue of (\ref{volumiuguali}),
(\ref{qfinitomaggioredizero}), (\ref{limitedebolevarphin}),
(\ref{uguale}), (\ref{1}), (\ref{convlambda}) and
(\ref{probleoverlinelambda}),  it results that
\begin{equation}\label{limann}\begin{array}{l}\displaystyle{
\lim_i a_{n_i}(v_{n_i},v_{n_i})= \overline\lambda[\overline
u,\overline u]_q -2\overline\lambda\sum_{k=1}^{\overline
k}[\overline u,u_{k}]_q^2+\sum_{k=1}^{\overline k}[\overline
u,u_{k}]_q^2\lambda_k=\overline\lambda.}
\end{array}\end{equation}

In what concerns the last term in (\ref{classicalpeoperties}),
equalities (\ref{deltahk}) provide that
\begin{equation}\nonumber\begin{array}{l}\displaystyle{
(v_n,v_n)_n=(\varphi_n-\sum_{k=1}^{\overline k} (\varphi_n,
u_{n,k})_n u_{n,k},\varphi_n-\sum_{j=1}^{\overline k} (\varphi_n,
u_{n,j})_n u_{n,j})_n= }\\\\
\displaystyle{(\varphi_n,\varphi_n)_n-2\sum_{k=1}^{\overline
k}(\varphi_n, u_{n,k})_n
(\varphi_n,u_{n,k})_n+\sum_{k,j=1}^{\overline k}(\varphi_n,
u_{n,k})_n (\varphi_n, u_{n,j})_n (u_{n,k},u_{n,j})_n=
}\\\\
\displaystyle{(\varphi_n,\varphi_n)_n-2\sum_{k=1}^{\overline
k}(\varphi_n, u_{n,k})_n^2+\sum_{k=1}^{\overline k}(\varphi_n,
u_{n,k})_n^2, = (\varphi_n,\varphi_n)_n-\sum_{k=1}^{\overline
k}(\varphi_n, u_{n,k})_n^2,\quad\forall n\in \mathbb{N}.}
\end{array}\end{equation}
Consequently, by virtue of  (\ref{volumiuguali}),
(\ref{qfinitomaggioredizero}), (\ref{limitedebolevarphin}),
(\ref{uguale}), (\ref{1}), (\ref{probleoverlinelambda}),  it
results that
\begin{equation}\label{lim()n}\begin{array}{l}\displaystyle{
\lim_i(v_{n_i},v_{n_i})_{n_i}=[\overline u,\overline u]_q-
\sum_{k=1}^{\overline k}[\overline u, u_{k}]_q^2=1.}
\end{array}\end{equation}

Finally, by passing to the limit in (\ref{classicalpeoperties}) as
$n_i\rightarrow+\infty$, convergences (\ref{limann}),
(\ref{convlambda}) and (\ref{lim()n}) entail that
\begin{equation}\nonumber
\overline\lambda\geq \lambda_{\overline k},
\end{equation}
which is in contradiction with (\ref{maggiore}). Hence, Problem
(\ref{probleoverlinelambda}) does not admit solution.

\noindent $\underline{\hbox {Step 3.}}$ Conclusion.

In Step 1, it is proved that $\{\lambda_{k}\}_{k\in
\mathbb{N}}\subset]0,+\infty[$ is an increasing and diverging
sequence of eigenvalues of Problem (\ref{problemalimite}),
$\{u_{k}\}_{k\in \mathbb{N}}$ is an orthonormal sequence in
$(H,[\cdot,\cdot]_q)$, $\{\lambda_k^{-\frac{1}{2}}u_{k}\}_{k\in
\mathbb{N}}$ is an orthonormal sequence in $(V,\alpha_q)$, and,
for every $k\in\mathbb{N}, $ $u_{k} $ is an eigenvector for
Problem (\ref{problemalimite}), with eigenvalue $\lambda_{k}$.

Indeed, the sequence $\{\lambda_{k}\}_{k\in \mathbb{N}}$ forms the
whole set of the eigenvalues of Problem (\ref{problemalimite}). In
fact, if $\overline\lambda\notin\{\lambda_{k}\}_{k\in \mathbb{N}}$
is another eigenvalue of Problem (\ref{problemalimite}) and
$\overline u\in V$ is a corresponding eigenvector (which can be
assumed $H$-normalized), it results that
$$\overline\lambda [\overline u,
u_k]_q=\alpha_q(\overline u,u_k)=\alpha_q(u_k,\overline
u)=\lambda_k[u_k,\overline u]_q,\quad\forall k\in\mathbb{N},$$
that is $$[\overline u, u_k]_q=0,\quad\forall k\in\mathbb{N},$$
which is in contradiction with the statement of Step 2.

It remains to prove that the set of the finite combinations of
elements of $\{\lambda_k^{-\frac{1}{2}}u_{k}\}_{k\in \mathbb{N}}$
is dense in $(V,\alpha_q)$, which provides that the set of the
finite combinations of elements of $\{u_{k}\}_{k\in \mathbb{N}}$
is dense in $(H,[\cdot,\cdot]_q)$, since $V$ is continuously
embedded into $H$, with dense inclusion.

It is well known that $(V,\alpha_q)$ has a Hilbert orthonormal
basis  of eigenvectors of Problem (\ref{limiteigenproblem}) (for
instance, see Th. 6.2-1 in \cite{RT}). Consequently,  by denoting
with $\{\mu_i\}_{i\in\mathbb{N}}$ the sequence of all different
eigenvalues of Problem (\ref{limiteigenproblem}) (i.e.
$\mu_i\neq\mu_j$ for $i\neq j$) and with
$\displaystyle{\{E_i\}_{i\in\mathbb{N}}}$ the sequence of the
 spaces of the associated eigenvectors, it results that the vectorial space generated by
$\displaystyle{\{E_i\}_{i\in\mathbb{N}}}$ is dense in
$(V,\alpha_q)$. Recall that, for every $i\in \mathbb{N}$, $E_i$
has finite dimension and $E_i$ is orthogonal to $ E_j$ (with
respect to $\alpha_q$) if $i\neq j$. Then, to conclude it is
enough to show that, for every $i\in \mathbb{N}$, an orthonormal
basis (with respect to $\alpha_q$) of $E_i$ is included in
$\{\lambda_k^{-\frac{1}{2}}u_{k}\}_{k\in \mathbb{N}}$. This last
property will be proved by arguing by contradiction. If it is not
true,  there exist $\overline i\in \mathbb{N}$ and $\overline u\in
E_{\overline i}$ such that $\alpha_q(\mu^{-\frac{1}{2}}_{\overline
i}\overline u,u_k)=0$ for every $k\in\mathbb{N}$, and
$\alpha_q(\mu^{-\frac{1}{2}}_{\overline i}\overline
u,\mu^{-\frac{1}{2}}_{\overline i}\overline u)=1$. That is
$\overline u$ is an eigenvector of Problem
(\ref{limiteigenproblem}), with eigenvalue $\mu_{\overline i}$,
such that $[\overline u,u_k]_q=0$ for every $k\in\mathbb{N}$, and
$[\overline u,\overline u]_q=1$, in contradiction with the
statement of Step 2.

In conclusion, since the sequence $\{\lambda_k\}_{k\in
\mathbb{N}}$ can be characterized by the min-max Principle (for
instance, see Th. 6.2-2 in \cite{RT}), for every $k\in
\mathbb{N}$, convergence (\ref{convlambda}) holds true for the
whole sequence $\{\lambda_{n,k}\}_{n\in \mathbb{N}}$.
\end{proof}

 \textit{Proof of Theorem \ref{TH1}.} As it is usual (see \cite{CD}), Problem
(\ref{opx}) can be reformulated
 on a fixed domain through an appropriate rescaling which
maps $\Omega_n$ into $\Omega=\omega \times ]-1,1[$.
 Namely,    by setting
\begin{equation}\nonumber
\widetilde{u}_{n,k}(x)=\left\{
\begin{array}{ll} \widetilde{u}^{a}_{n,k}(x',x_3)=
 U_{n,k}(r_nx',x_3),\quad (x',x_3)\hbox{ a.e. in
}\Omega^a=\omega\times ]0,1[, \\\\
 \widetilde{u}^{b}_{n,k}(x',x_3)= U_{n,k}(x',h_n x_3), \quad
(x',x_3)\hbox{ a.e. in }\Omega^b=\omega\times ]-1,0[,
\end{array}\right.\forall n ,k\in \mathbb N,
\end{equation}
it results that, for every  $n \in \mathbb N$,
$\{\lambda_{n,k}\}_{k \in \mathbb N}$ forms the set of all  the
eigenvalues of Problem (\ref{eigenproblem}),
$u_{n,k}={r_n^\frac{N-1}{2}}\widetilde{u}_{n,k}\in V_n$ is an
eigenvector of (\ref{eigenproblem}) with eigenvalue
$\lambda_{n,k}$, $\{u_{n,k}\}_{k \in \mathbb N}$ is a
$H_n$-Hilbert orthonormal basis  and
$\left\{\lambda_{n,k}^{-\frac{1}{2}}u_{n,k}\right\}_{k \in \mathbb
N}$ is a $V_n$-Hilbert orthonormal basis, where $V_n$ and $H_n$
are defined at the beginning of this section. Then, Theorem
\ref{TH1} is an immediate consequence  of Theorem
\ref{maintheorem}.$\hfill\Box$

\section{The case $\displaystyle{\lim_n\frac{h_n}{r_n^{N-1}}=0}$\label{Sec3}}

For every $n ,k\in \mathbb N$, let $\lambda_{n,k}$ be as in
Section \ref{mainresult}. The aim of this  section is  to
investigate the limit behavior, as $n\rightarrow+\infty$, of
$\{\lambda_{n,k}\}_{n \in \mathbb N}$,  under the assumption:
\begin{equation}\label{q=0}\lim_n\frac{h_n}{r_n^{N-1}}=0.\end{equation}

Let $H_0$ be the space $H$ given in (\ref{H}) equipped with the
 inner product $ [\cdot,\cdot]_{1}$ defined
 by (\ref{prodottoscalarelim}) with  $q=1$. Moreover,
let
\begin{equation}\label{Vb0}V_{0b}=\left\{\begin{array}{ll}\begin{array}{l}
\Big\{ v^b \in  H^1 (\Omega^b):
 v^b \hbox{ is independent of
}x_N,\, v^b=0 \hbox{ on } \partial\omega,\\
v^b(0')=0\Big\},\end{array}&\hbox{ if }  N=2; \\\\\\
\begin{array}{l}\Big\{ v^b \in  H^1 (\Omega^b):
 v^b \hbox{ is independent of
}x_N,\, v^b=0 \hbox{ on } \partial\omega\Big\},\end{array} &\hbox{
if } 3\leq N
\end{array}\right.\end{equation}
equipped with the  $H_0^1(\omega)$-norm, and $V_0$ be the space
\begin{equation}\label{V0}V_0=\left\{v^a\in H^1(\Omega^a):v^a \hbox{ is independent of
}x',\,v^a(1)=0\right\}\times V_{0b}\end{equation} equipped with
the  inner product $\alpha_1$ defined by (\ref{formabilineare})
with $q=1$. Remark that  $V_0$ is continuously and compactly
embedded into $H_0$, with dense inclusion.

Let
\begin{equation}\label{defV0b}\widetilde{V}_{0b}=\Big\{ v^b \in  W_0^{1,\infty} (\omega):
v^b=0\hbox{ in a  neighbourhood of }0'\Big\},\end{equation} then
the following density result holds true:
\begin{Proposition}\label{densità}
$\widetilde{V}_{0b}$ is dense in ${V}_{0b}.$
\end{Proposition}
\begin{proof}[Proof]We sketch the proof.

The density result is obvious if $N=2$.

Assume $N> 2$ and let $\{\varepsilon_n\}_{n\in \mathbb{N}}\subset
\mathbb{R}$, $\{\eta_n\}_{n\in \mathbb{N}}\subset \mathbb{R}$ and
$\{\varphi_n\}_{n\in \mathbb{N}}\subset C(\mathbb{R}^{N-1})$ three
sequences such that, for every $ n\in \mathbb{N}$,
$0<\varepsilon_n<\eta_n<\hbox{dist}(0',\partial\omega)$,
$\varphi_n=1$ in $\{x'\in \mathbb{R}^{N-1}: \vert x'\vert\leq
\varepsilon_n\}$, $\varphi_n=0$ in $\mathbb{R}^{N-1}-\{x'\in
\mathbb{R}^{N-1}: \vert x'\vert< \eta_n\}$, $\{\varphi_n\}_{n\in
\mathbb{N}}\subset C^1(\{x'\in \mathbb{R}^{N-1}: \varepsilon_n\leq
\vert x'\vert\leq \eta_n\})$, $0\leq \varphi_n\leq 1$, and
$\displaystyle{\lim_n\eta_n=0=\lim_n\int_{\{x'\in
\mathbb{R}^{N-1}: \varepsilon_n<\vert x'\vert< \eta_n\}}\vert D
\varphi_n\vert^2dx'}$ (for the existence of such  sequences, see
 Proposition  3.1 in \cite{GGLM1}).

 Since $C^\infty_0(\omega)$ is dense in ${V}_{0b}$, it is enough to prove that $\widetilde{V}_{0b}$ is
dense in $C^\infty_0(\omega)$ with respect to the
$H_0^1(\omega)$-norm. For $v\in C^\infty_0(\omega)$, set
$v_n=(1-\varphi_n)v\in \widetilde{V}_{0b}$ for every $n\in
\mathbb{N}$. Since $v=\varphi_nv+(1-\varphi_n)v$  for every $n\in
\mathbb{N}$, and
\begin{equation}\nonumber\begin{array}{l}\displaystyle{\lim_n\int_{\omega} \vert D
(v\varphi_n)\vert^2dx'\leq}\\\\\displaystyle{ \Vert
v\Vert^2_{W^{1,\infty}(\omega)}\lim_n \left(\hbox{meas} \{x'\in
\mathbb{R}^{N-1}: \vert x'\vert< \eta_n\}+\int_{\{x'\in
\mathbb{R}^{N-1}: \varepsilon_n<\vert x'\vert< \eta_n\}}\vert D
\varphi_n\vert^2dx'\right)=0,}\end{array}\end{equation} it follows
that $v_n\rightarrow v$ strongly in $ H^1_0(\omega)$.
\end{proof}

Consider the eigenvalue problem:
\begin{equation}\label{limiteigenproblem0}\left\{\begin{array}{l}
\displaystyle{u\in  V_0,}\\\\
\displaystyle{\alpha_1(u, v)=\lambda[u, v]_{1},\quad \forall v\in
{ V_0}},
\end{array}\right.\end{equation}
then, if  (\ref{q=0}) holds true, the following convergence result
yields:
\begin{Theorem}\label{theorem0}
For every $n \in \mathbb N$, let $\{\lambda_{n,k}\}_{k \in \mathbb
N}$ be an increasing diverging sequence of  all the eigenvalues of
Problem (\ref{prodottoscalaren})$\div$(\ref{eigenproblem}), and
$\{\overline u_{n,k}\}_{k \in \mathbb N}$ be a Hilbert orthonormal
basis for the space $H_n$ equipped with the inner product:
$\displaystyle{\frac{r_n^{N-1}}{h_n}(\cdot, \cdot)_n}$, such that
$\left\{\lambda_{n,k}^{-\frac{1}{2}}\overline u_{n,k}\right\}_{k
\in \mathbb N}$ is a Hilbert orthonormal basis for the space $V_n$
equipped with the inner product
$\displaystyle{\frac{r_n^{N-1}}{h_n}a_n(\cdot, \cdot)}$ and, for
every $k\in \mathbb N$, $\overline u_{n,k}$ is an eigenvector of
Problem (\ref{prodottoscalaren})$\div$(\ref{eigenproblem}) with
eigenvalue $\lambda_{n,k}$.   Assume (\ref{hrzero}) and
(\ref{q=0}).

Then, there exists an increasing diverging sequence of positive
numbers $\{\lambda_{k}\}_{k \in \mathbb N}$ such that
$$\lim_n\lambda_{{n},k}=\lambda_k, \quad\forall k \in \mathbb N,$$
and $\{\lambda_{k}\}_{k \in \mathbb N}$ is the set of all  the
eigenvalues of Problem (\ref{Vb0}), (\ref{V0}),
(\ref{limiteigenproblem0}).
 Moreover,
there exists an increasing sequence of positive integer numbers
$\{n_i\}_{i \in \mathbb N}$ and an $H_0$-Hilbert orthonormal basis
$\{u_{k}=(u_{k}^a,u_{k}^b)\}_{k \in \mathbb N}$ (depending
possibly on the selected subsequence $\{n_i\}_{i \in \mathbb N}$)
such that, for every $k\in \mathbb N$, $u_{k}\in V_0$ is an
eigenvector of Problem (\ref{Vb0}), (\ref{V0}),
(\ref{limiteigenproblem0}), with eigenvalue $\lambda_{k}$, and
\begin{equation}\label{00}\displaystyle{
\frac{r_{n_i}^{\frac{N-1}{2}}}{h_{n_i}^{\frac{1}{2}}}\overline
u_{n_i,k}^a\rightarrow u_k^a\hbox{ strongly in
}H^1(\Omega^a),\quad \overline u_{n_i,k}^b\rightarrow u_k^b,\hbox{
strongly in }H^1(\Omega^b),\quad\forall k \in \mathbb N,}
\end{equation}
as $i\rightarrow+\infty$,
\begin{equation}\label{limitedeboleautovettoribiszcx}\left\{\begin{array}{l}\displaystyle{
\frac{1}{r_{n}}\frac{r_{n}^{\frac{N-1}{2}}}{h_{n}^{\frac{1}{2}}}D_{x'}\overline
u_{n,k}^a\rightarrow 0'\hbox{ strongly in
}(L^2(\Omega^a))^{N-1},}\\\\ \displaystyle{
\frac{1}{h_{n}}\partial_{x_N}\overline u_{n,k}^b\rightarrow
0\hbox{ strongly in }L^2(\Omega^b),}
\end{array}\right. \quad\forall k \in \mathbb N,\end{equation}
as $n\rightarrow+\infty$. Furthermore,
$\{\lambda_k^{-\frac{1}{2}}u_{k}\}_{k\in \mathbb{N}}$ is a
$V_0$-Hilbert orthonormal basis.
\end{Theorem}
\begin{proof}[Proof]We sketch the proof.

By multiplying  (\ref{eigenproblem}) by
$\displaystyle{\frac{r_n^{N-1}}{h_n}}$, it results that
\begin{equation}\label{eigenproblemperun,kbis}\left\{\begin{array}{l}
\displaystyle{\overline u_{n,k}\in { V}_n,}\\\\
\displaystyle{\frac{r_n^{N-1}}{h_n}a_n(\overline u_{n,k},
v)=\lambda_{n,k}\frac{r_n^{N-1}}{h_n}(\overline u_{n,k},
v)_n,\quad \forall v\in { V}_n,}
\end{array}\right.\quad\forall n,k \in \mathbb
N.\end{equation} Moreover, since $\{\overline u_{n,k}\}_{k \in
\mathbb N}$ is a Hilbert orthonormal basis for the space $H_n$
equipped with the inner product
$\displaystyle{\frac{r_n^{N-1}}{h_n}(\cdot, \cdot)_n}$, one has
that
\begin{equation}\label{deltahkbis}
\frac{r_n^{N-1}}{h_n}(\overline u_{n,k}, \overline
u_{n,h})_n=\delta_{h,k},\quad \forall n,k,h, \in \mathbb N,
\end{equation}
where $\delta_{h,k}$ is defined in (\ref{delta})

By choosing $v=\overline u_{n,k}$ in
(\ref{eigenproblemperun,kbis}), and by taking into account
(\ref{deltahkbis}), Proposition \ref{stimeautovalori} and Remark
\ref{remarkstimaautov}, it follows that
\begin{equation}\nonumber
\frac{r_n^{N-1}}{h_n}a_n(\overline u_{n,k}, \overline
u_{n,k})=\lambda_{n,k}\frac{r_n^{N-1}}{h_n}(\overline u_{n,k},
\overline u_{n,k})_n=\lambda_{n,k}\leq 2^kk^2\pi^2,\quad\forall n
,k \in \mathbb N.
\end{equation}
Consequently, by virtue of  definition (\ref{formabilinearen}) of
$a_n$ and of  assumption (\ref{hrzero}), and  by using  a diagonal
argument, it is easily seen that there exists an increasing
sequence of positive integer numbers $\{n_i\}_{i \in \mathbb N}$,
an increasing sequence of no negative numbers $\{\lambda_{k}\}_{k
\in \mathbb N}$, a sequence
$\left\{u_{k}=(u_{k}^a,u_{k}^b)\right\}_{k \in \mathbb N}\subset
H^1(\Omega^a)\times H^1(\Omega^b)$ (depending possibly on the
selected subsequence $\{n_i\}_{i \in \mathbb N}$), with $u_{k}^a$
independent of $x'$, $u_{k}^b$ independent of $x_N$,
$u_{k}^a(1)=0$, $u_{k}^b=0$ on $\partial\omega$,  and a sequence
$\left\{(\xi_{k}^a,\xi_{k}^b)\right\}_{k \in \mathbb N}\subset
(L^2(\Omega^a))^{N-1}\times L^2(\Omega^b)$ (depending possibly on
the selected subsequence $\{n_i\}_{i \in \mathbb N}$) such that
\begin{equation}\label{convlambdabis}\lim_i\lambda_{{n_i},k}=\lambda_k, \quad\forall k \in \mathbb N,\end{equation}
\begin{equation}\label{limitedeboleautovettoribis}\left\{\begin{array}{l}\displaystyle{
\frac{r_{n_i}^{\frac{N-1}{2}}}{h_{n_i}^{\frac{1}{2}}}\overline
u_{n_i,k}^a\rightharpoonup u_k^a\hbox{ weakly in
}H^1(\Omega^a)\hbox{ and strongly in
}L^2(\Omega^a),}\\\\
\displaystyle{ \overline u_{n_i,k}^b\rightharpoonup u_k^b\hbox{
weakly in }H^1(\Omega^b)\hbox{ and strongly in }L^2(\Omega^b),}
\end{array}\right. \quad\forall k \in \mathbb N,\end{equation}
\begin{equation}\label{limitedeboleautovettoribisz}\left\{\begin{array}{l}\displaystyle{
\frac{r_{n_i}^{\frac{N-1}{2}}}{h_{n_i}^{\frac{1}{2}}}\frac{1}{r_{n_i}}D_{x'}\overline
u_{n_i,k}^a\rightharpoonup \xi^a_k\hbox{ weakly in
}(L^2(\Omega^a))^{N-1},}\\\\
\displaystyle{ \frac{1}{h_{n_i}}\partial_{x_N}\overline
u_{n_i,k}^b\rightharpoonup  \xi^b_k\hbox{ weakly in
}L^2(\Omega^b),}
\end{array}\right. \quad\forall k \in \mathbb N,\end{equation}
as $i\rightarrow +\infty$.

To obtain $u^k\in V_0$, for every $k\in \mathbb{N}$, it remains to
prove that
\begin{equation}\label{ubk=0}u_{k}^b(0')=0, \quad\forall k\in \mathbb{N},\quad \hbox { if }N=2.\end{equation} Indeed,
statement (2.4) in \cite{GGLM1} (which
holds true also if $q=0$!) gets
\begin{equation}\label{N=2}
\lim_i\int_\omega \overline
u_{n_i,k}^b(r_{n_i}x',0)dx'=\vert\omega\vert
u_{k}^b(0'),\quad\forall k \in \mathbb N, \quad\hbox{ if }N=2.
\end{equation}
Then, by combining  (\ref{N=2}) with the first line in
(\ref{limitedeboleautovettoribis}), and by taking into account
that $u_{n_i,k}^a(x',0)= u_{n_i,k}^b(r_{n_i}x',0)$ for $x'$ a.e.
in $\omega$ and  (\ref{q=0}), one obtains that
\begin{equation}\nonumber\begin{array}{l}\displaystyle{\vert\omega\vert
u_{k}^b(0')=\lim_i\int_\omega \overline
u_{n_i,k}^b(r_{n_i}x',0)dx'=\lim_i\int_\omega \overline
u_{n_i,k}^a(x',0)dx'=}\\\\\displaystyle{\lim_i\left(\frac{h_{n_i}^{\frac{1}{2}}}{r_{n_i}^{\frac{N-1}{2}}}\int_\omega
\frac{r_{n_i}^{\frac{N-1}{2}}}{h_{n_i}^{\frac{1}{2}}}\overline
u_{n_i,k}^a(x',0)dx'\right)=} 0\vert\omega\vert
u^a_k(0)=0,\quad\forall k \in \mathbb N, \quad\hbox{ if
}N=2,\end{array} \end{equation} that is (\ref{ubk=0}).

Now, by passing to the limit, as $n_i\rightarrow +\infty$, in
(\ref{eigenproblemperun,kbis}), with $v=\left\{\begin{array}{ll}0&
\hbox{ in }\Omega^a,\\v^b &\hbox{ in }\Omega^b\end{array}\right.$
and $v^b\in \widetilde{V}_{0b}$, and by making use of
(\ref{convlambdabis}) and of the second line of
(\ref{limitedeboleautovettoribis}), it turns out that
\begin{equation}\nonumber
\int_{\omega}D_{x'}u^b_kD_{x'}v^bdx'=\lambda_k\int_{\omega}u^b_kv^bdx',\quad\forall
v^b\in \widetilde{V}_{0b},\quad\forall k \in \mathbb N,
\end{equation}
and, consequently, by virtue of Proposition  \ref{densità},
\begin{equation}\label{equainOb}
\int_{\omega}D_{x'}u^b_kD_{x'}v^bdx'=\lambda_k\int_{\omega}u^b_kv^bdx',\quad\forall
v^b\in {V}_{0b},\quad\forall k \in \mathbb N.
\end{equation}

On the other hand, by passing to the limit, as $n_i\rightarrow
+\infty$, in
\begin{equation}\nonumber\displaystyle{\frac{r_n^{\frac{N-1}{2}}}{h_n^{\frac{1}{2}}}a_n(\overline
u_{n,k},
v)=\lambda_{n,k}\frac{r_n^{\frac{N-1}{2}}}{h_n^{\frac{1}{2}}}(\overline
u_{n,k}, v)_n,\quad \forall v\in { V}_n,\quad\forall k \in \mathbb
N,}\end{equation} with $v=\left\{\begin{array}{ll}v^a & \hbox{ in
}\Omega^a,\\v^b &\hbox{ in }\Omega^b\end{array}\right.$ and
$v^a\in H^1(]0,1[)$, $v^a(1)=0$, $v^b\in C^1_0(\omega)$,
$v^b=v^a(0)$ in a neighbourhood of $0'$, and by making use of
(\ref{convlambdabis}) and of (\ref{limitedeboleautovettoribis}),
it turns out that
\begin{equation}\label{equainOa}
\vert\omega\vert\int_0^1\partial_{x_N}u^a_k\partial_{x_N}v^adx_N=\lambda_k\vert\omega\vert\int_0^1u^a_kv^adx_N,\,\,\forall
v^a\in H^1(]0,1[):v^a(1)=0, \,\,\forall k \in \mathbb N.
\end{equation}

By adding (\ref{equainOb}) to (\ref{equainOa}), one obtains that
\begin{equation}\label{problemalimitexyz}\left\{\begin{array}{l}
\displaystyle{u_{k}\in { V_0},}\\\\
\displaystyle{\alpha_1(u_{k}, v)=\lambda_{k}[u_{k}, v]_1,\quad
\forall v\in { V_0},}
\end{array}\right.\quad\forall k\in \mathbb N.\end{equation}

Finally, by passing to the limit, as $n_i\rightarrow +\infty$, in
(\ref{deltahkbis}), and by making use of
(\ref{limitedeboleautovettoribis}) and (\ref{problemalimitexyz}),
it turns out that
\begin{equation}\label{problemalimitexyzç}
\alpha_1(\lambda_{k}^{-\frac{1}{2}}u_{k},\lambda_{h}^{-\frac{1}{2}}
u_{h})=[u_{k}, u_{h}]_1=\delta_{h,k},\quad\forall k,h\in \mathbb
N.\end{equation} Moreover,  by arguing as in  the proof of Theorem
\ref{maintheorem}, one proves that
$$\displaystyle{\lim_{k}\lambda_{k}=+\infty.}$$

Now, let us identify  $\xi^a_k$ and $\xi^b_k$. By choosing
$v=\overline u_{n,k}$ as test function in
(\ref{eigenproblemperun,kbis}), by replacing $n$ with $n_i$, by
passing to the $\liminf$  as $i$ diverges, by taking into account
of (\ref{convlambdabis}), (\ref{limitedeboleautovettoribis}),
(\ref{limitedeboleautovettoribisz}) and
(\ref{problemalimitexyzç}), and by using a l.s.c. argument it
results that
\begin{equation}\nonumber\int_{\Omega^a}\vert\xi^a_k\vert^2+\vert
\partial_{x_N}u^a_k\vert^2dx+\int_{\Omega^b}
\vert D_{x'}u^b_k\vert^2+\vert\xi^b_k\vert^2dx\leq
\int_{\Omega^a}\vert
\partial_{x_N}u^a_k\vert^2dx+\int_{\Omega^b}
\vert D_{x'}u^b_k\vert^2dx,\,\,\forall k\in \mathbb N,
\end{equation}
which provides that
\begin{equation}\label{xa=0}\xi^a_k=0,\quad\xi^b_k=0,\quad\forall k\in \mathbb
N.
\end{equation}

The  strong convergences in (\ref{00}) and
(\ref{limitedeboleautovettoribiszcx}) follows from
(\ref{limitedeboleautovettoribis}),
(\ref{limitedeboleautovettoribisz}), (\ref{xa=0}) and from the
convergence of the energies:
\begin{equation}\lim_i
\frac{r_{n_i}^{N-1}}{h_{n_i}}a_{n_i}(\overline u_{{n_i},k},
\overline
u_{{n_i},k})=\lim_i\lambda_{n_i,k}=\lambda_k=\alpha_1(u_{k},
u_{k}),\quad \forall k \in \mathbb N.
\end{equation}

 As in Step 2 of the proof of Theorem \ref{maintheorem}, by arguing by
contradiction, one can show that   there not exist $(\overline u,
\overline \lambda)\in V_0\times\mathbb{R}$ satisfying the
following problem:
\begin{equation}\label{probleoverlinelambdakj}\left\{\begin{array}{l}
\displaystyle{\overline
u\in { V_0},}\\\\
\displaystyle{\alpha_1(\overline u, v)=\overline\lambda [\overline
u, v]_1,\quad
\forall v\in { V_0},}\\\\

[\overline u, u_k]_1=0,\quad\forall k\in\mathbb N\\\\

[\overline u,\overline u ]_1=1.
\end{array}\right.\end{equation}
Precisely, assume that there exists $(\overline u, \overline
\lambda)\in V_0\times\mathbb{R}$ satisfying
(\ref{probleoverlinelambdakj}). Let
 ${\overline k}\in \mathbb{N}$ be such that
\begin{equation}\nonumber
\overline\lambda< \lambda_{\overline k},
\end{equation} and,
for every $n\in \mathbb{N}$, let $\varphi_n$ be the solution of
the following problem:
\begin{equation}\nonumber\left\{\begin{array}{l}
\displaystyle{\varphi_n\in { V}_n,}\\\\
\displaystyle{a_n(\varphi_n,
v)=\overline\lambda\left(\left(\frac{h_n^{\frac{1}{2}}}{r_n^{\frac{N-1}{2}}}\overline
u^a,\overline u^b\right), v\right)_n,\quad \forall v\in { V}_n,}
\end{array}\right.\end{equation}
where $\overline{u}=(\overline u^a,\overline u^b)$. Then, it is
easy to prove that
\begin{equation}\nonumber\displaystyle{
\frac{r_{n}^{\frac{N-1}{2}}}{h_{n}^{\frac{1}{2}}}\varphi_{n}^a\rightharpoonup
\overline {u}^a\hbox{ weakly in }H^1(\Omega^a),\quad
\varphi_{n}^b\rightharpoonup \overline{u}^b\hbox{ weakly in
}H^1(\Omega^b),}
\end{equation}
as $n\rightarrow +\infty$. Moreover, for every $n\in \mathbb{N}$,
set
$$v_n=\varphi_n-\sum_{i=1}^{\overline k}
\frac{r_n^{N-1}}{h_n}(\varphi_n, u_{n,i})_n u_{n,i}\in V_n.$$
Then, by proceeding as in the proof of Step 2 of Theorem
\ref{maintheorem}, but  by taking care to replace $(\cdot,
\cdot)_n$ with $\displaystyle{\frac{r_n^{N-1}}{h_n}(\cdot,
\cdot)_n}$ and $a_n(\cdot, \cdot)$ with
$\displaystyle{\frac{r_n^{N-1}}{h_n}a_n(\cdot, \cdot)}$, one
reaches a contradiction.

To complete the proof, one can argue as  Step 3 of the proof of
Theorem \ref{maintheorem}.\end{proof}

\begin{Remark} \label{remm}For every $n \in \mathbb N$, let $\{\lambda_{n,k}\}_{k \in \mathbb
N}$ be an increasing diverging sequence of  all the eigenvalues of
Problem (\ref{prodottoscalaren})$\div$(\ref{eigenproblem}), and
$\{ u_{n,k}\}_{k \in \mathbb N}$ be a Hilbert orthonormal basis
for the space $H_n$, such that
$\left\{\lambda_{n,k}^{-\frac{1}{2}} u_{n,k}\right\}_{k \in
\mathbb N}$ is a Hilbert orthonormal basis for the space $V_n$,
and, for every $k\in \mathbb N$, $u_{n,k}$ is an eigenvector of
(\ref{eigenproblem}) with eigenvalue $\lambda_{n,k}$.

By setting $\displaystyle{\overline
u_{n,k}=\frac{h_n^{\frac{1}{2}}}{r_n^{\frac{N-1}{2}}}u_{n,k}}$,
 it turns out that, for every $n\in \mathbb{N}$, $\{\overline
u_{n,k}\}_{k \in \mathbb N}$ is a Hilbert orthonormal basis for
the space $H_n$ equipped with the inner product
$\displaystyle{\frac{r_n^{N-1}}{h_n}(\cdot, \cdot)_n}$,
$\left\{\lambda_{n,k}^{-\frac{1}{2}}\overline u_{n,k}\right\}_{k
\in \mathbb N}$ is a Hilbert orthonormal basis for the space $V_n$
equipped with the inner product
$\displaystyle{\frac{r_n^{N-1}}{h_n}a_n(\cdot, \cdot)}$, and, for
every $k\in \mathbb N$, $\overline u_{n,k}$ is an eigenvector of
(\ref{eigenproblem}) with eigenvalue $\lambda_{n,k}$.

Then,  by applying Theorem \ref{theorem0}, it follows that
$$\lim_n\lambda_{{n},k}=\lambda_k, \quad\forall k \in \mathbb N,$$
\begin{equation}\nonumber\displaystyle{
u_{n_i,k}^a\rightarrow u_k^a\hbox{ strongly in
}H^1(\Omega^a),\quad
\frac{h_{n_i}^{\frac{1}{2}}}{r_{n_i}^{\frac{N-1}{2}}}
u_{n_i,k}^b\rightarrow u_k^b,\hbox{ strongly in
}H^1(\Omega^b),\quad\forall k \in \mathbb N,}
\end{equation}
as $i\rightarrow+\infty$,
\begin{equation}\nonumber\displaystyle{\frac{1}{r_n}D_{x'}u_{n,k}^a\rightarrow 0'\hbox{ strongly in
}(L^2(\Omega^a))^{N-1},\,\,
\frac{1}{h_n^{\frac{1}{2}}r_n^{\frac{N-1}{2}}}
\partial_{x_N}u_{n,k}^b\rightarrow 0,\hbox{ strongly in
}L^2(\Omega^b),\,\,\forall k \in \mathbb N,}
\end{equation}
as $n\rightarrow+\infty$, where $\{\lambda_{k}\}_{k \in \mathbb
N}$ and $\{(u_k^a,u_k^b)\}_{k \in \mathbb N}$ are given by Theorem
\ref{theorem0}.
\end{Remark}

\textit{Proof of Theorem \ref{TH2}.}  Theorem \ref{TH2} follows
immediately from Remark \ref{remm}, by arguing as in the proof of
Theorem \ref{TH1} at the end of Section
\ref{mainresult}.$\hfill\Box$

\section{The case $\displaystyle{\lim_n\frac{h_n}{r_n^{N-1}}=+\infty}$\label{Sec4}}

For every $n ,k\in \mathbb N$, let $\lambda_{n,k}$ be as in
Section \ref{mainresult}. The aim of this  section is  to
investigate the asymptotic behavior, as $n\rightarrow+\infty$, of
$\{\lambda_{n,k}\}_{n \in \mathbb N}$,  under the assumption:
\begin{equation}\label{q=+infty}\lim_n\frac{h_n}{r_n^{N-1}}=+\infty.\end{equation}
Indeed, in this case,  only a partial result is obtained.
Precisely, for $N=2$, the limit eigenvalue problem is completely
derived; while, for $N\geq 3$, it is obtained only under the
additional assumption:
\begin{equation}\label{ulterioreipotesi}
\left\{\begin{array}{l}h_n\ll -r_n^2\log r_n, \hbox{ if } N=3,\\\\
h_n\ll r_n^2,\hbox{ if } N\geq 4 .\end{array}\right.\end{equation}

As in Section \ref{Sec3}, let $H_0$ be the space $H$ given in
(\ref{H}) equipped with the
 inner product $ [\cdot,\cdot]_{1}$ defined
 by (\ref{prodottoscalarelim}) with  $q=1$. Moreover,
let $V_\infty$ be the  space:
\begin{equation}\label{limitspaceinf}V_\infty=\left\{\begin{array}{ll}\begin{array}{l}
\Big\{v=(v^a, v^b) \in H^1 (\Omega^a)\times H^1 (\Omega^b):\\\\
v^a \hbox{ is independent of }x',\, v^b \hbox{ is independent of
}x_N,\\\\ v^a(1)=0,\,\, v^a(0)=0, \\\\
v^b=0 \hbox{ on } \partial\omega\Big\},\end{array}&\quad\quad\hbox{ if }  N=2; \\\\\\
\begin{array}{l}\Big\{v=(v^a, v^b) \in H^1 (\Omega^a)\times H^1
(\Omega^b):\\\\ v^a \hbox{ is
independent of }x',\, v^b \hbox{ is independent of }x_N,\\\\
v^a(1)=0,\,\, v^b=0 \hbox{ on } \partial\omega \Big\},\end{array}
&\quad\quad\hbox{ if } 3\leq N;
\end{array}\right.\end{equation}
equipped with the  inner product $\alpha_1$ defined by
(\ref{formabilineare}) with $q=1$. Remark that  $V_\infty$ is
continuously and compactly embedded into $H_0$, with dense
inclusion.

Consider the
 following eigenvalue problem:
 \begin{equation}\label{limiteigenprobleminf}\left\{\begin{array}{l}
\displaystyle{u\in  V_\infty,}\\\\
\displaystyle{\alpha_1(u, v)=\lambda[u, v]_{1},\quad \forall v\in
{ V_\infty}},
\end{array}\right.\end{equation}
then, if  (\ref{q=+infty}) and (\ref{ulterioreipotesi}) hold true,
the following convergence result yields:
\begin{Theorem}\label{theorem+infty}
For every $n \in \mathbb N$, let $\{\lambda_{n,k}\}_{k \in \mathbb
N}$ be an increasing diverging sequence of  all the eigenvalues of
Problem (\ref{prodottoscalaren})$\div$(\ref{eigenproblem}), and
$\{ u_{n,k}\}_{k \in \mathbb N}$ be a $H_n$-Hilbert orthonormal
basis  such that $\left\{\lambda_{n,k}^{-\frac{1}{2}}
u_{n,k}\right\}_{k \in \mathbb N}$ is a $V_n$-Hilbert orthonormal
basis  and, for every $k\in \mathbb N$, $ u_{n,k}$ is an
eigenvector of Problem
(\ref{prodottoscalaren})$\div$(\ref{eigenproblem}) with eigenvalue
$\lambda_{n,k}$. Assume (\ref{hrzero}), (\ref{q=+infty}) and, for
$N\geq 3$, also (\ref{ulterioreipotesi}).

Then, there exists an increasing diverging sequence of positive
numbers $\{\lambda_{k}\}_{k \in \mathbb N}$ such that
$$\lim_n\lambda_{{n},k}=\lambda_k, \quad\forall k \in \mathbb N,$$
and $\{\lambda_{k}\}_{k \in \mathbb N}$ is the set of all  the
eigenvalues of Problem (\ref{limitspaceinf}),
(\ref{limiteigenprobleminf}).
 Moreover,
there exists an increasing sequence of positive integer numbers
$\{n_i\}_{i \in \mathbb N}$ and an $H_0$-Hilbert orthonormal basis
$\{u_{k}=(u_{k}^a,u_{k}^b)\}_{k \in \mathbb N}$ (depending
possibly on the selected subsequence $\{n_i\}_{i \in \mathbb N}$)
such that, for every $k\in \mathbb N$, $u_{k}\in V_\infty$ is an
eigenvector of Problem (\ref{limitspaceinf}),
(\ref{limiteigenprobleminf}),
 with eigenvalue
$\lambda_{k}$, and
\begin{equation}\label{abinf}\displaystyle{
 u_{n_i,k}^a\rightarrow u_k^a\hbox{ strongly in
}H^1(\Omega^a),\quad
\frac{h_{n_i}^{\frac{1}{2}}}{r_{n_i}^{\frac{N-1}{2}}}
u_{n_i,k}^b\rightarrow u_k^b,\hbox{ strongly in
}H^1(\Omega^b),\quad\forall k \in \mathbb N,}
\end{equation}
as $i\rightarrow+\infty$.
\begin{equation}\label{limitedeboleautovettoribiszcxinf}\left\{\begin{array}{l}\displaystyle{
\frac{1}{r_{n}}D_{x'}u_{n,k}^a\rightarrow 0'\hbox{ strongly in
}(L^2(\Omega^a))^{N-1},}\\\\ \displaystyle{
\frac{1}{h_{n}}\frac{h_{n}^{\frac{1}{2}}}{r_{n}^{\frac{N-1}{2}}}\partial_{x_N}
u_{n,k}^b\rightarrow 0\hbox{ strongly in }L^2(\Omega^b),}
\end{array}\right. \quad\forall k \in \mathbb N,\end{equation}
as $n\rightarrow+\infty$. Furthermore,
$\{\lambda_k^{-\frac{1}{2}}u_{k}\}_{k\in \mathbb{N}}$ is a
$V_\infty$-Hilbert orthonormal basis.
\end{Theorem}

\begin{proof}[Proof]We sketch the proof.

It results that
\begin{equation}\label{eigenproblemperun,k+infty}\left\{\begin{array}{l}
\displaystyle{u_{n,k}\in { V}_n,}\\\\
\displaystyle{a_n(u_{n,k}, v)=\lambda_{n,k}(u_{n,k}, v)_n,\quad
\forall v\in { V}_n,}
\end{array}\right.\quad\forall n,k \in \mathbb
N,\end{equation}
\begin{equation}\label{deltahk+infty}
(u_{n,k}, u_{n,h})_n=\delta_{h,k},\quad \forall n,k,h, \in \mathbb
N,
\end{equation}
where $\delta_{h,k}$ is defined in (\ref{delta}).

By choosing $v=u_{n,k}$ in (\ref{eigenproblemperun,k+infty}), and
by taking into account (\ref{deltahk+infty}) and Proposition
\ref{stimeautovalori}, it follows that
\begin{equation}\label{aa}
a_n(u_{n,k}, u_{n,k})=\lambda_{n,k}(u_{n,k},
u_{n,k})_n=\lambda_{n,k}\leq 2^kk^2\pi^2,\quad\forall n ,k \in
\mathbb N.
\end{equation}
Consequently, by virtue of  definition (\ref{formabilinearen}) of
$a_n$ and of  assumption (\ref{hrzero}), and  by using  a diagonal
argument, it is easily seen that there exists an increasing
sequence of positive integer numbers $\{n_i\}_{i \in \mathbb N}$,
an increasing sequence of no negative numbers $\{\lambda_{k}\}_{k
\in \mathbb N}$, a sequence
$\left\{u_{k}=(u_{k}^a,u_{k}^b)\right\}_{k \in \mathbb N}\subset
H^1(\Omega^a)\times H^1(\Omega^b)$ (depending possibly on the
selected subsequence $\{n_i\}_{i \in \mathbb N}$), with $u_{k}^a$
independent of $x'$, $u_{k}^b$ independent of $x_N$,
$u_{k}^a(1)=0$, $u_{k}^b=0$ on $\partial\omega$, and a sequence
$\left\{(\xi_{k}^a,\xi_{k}^b)\right\}_{k \in \mathbb N}\subset
(L^2(\Omega^a))^{N-1}\times L^2(\Omega^b)$ (depending possibly on
the selected subsequence $\{n_i\}_{i \in \mathbb N}$) such that
\begin{equation}\label{convlambda+infty}\lim_i\lambda_{{n_i},k}=\lambda_k, \quad\forall k \in \mathbb N,\end{equation}
\begin{equation}\label{limitedeboleautovettori+infty}\left\{\begin{array}{l}\displaystyle{
 u_{n_i,k}^a\rightharpoonup u_k^a\hbox{ weakly in
}H^1(\Omega^a)\hbox{ and strongly in
}L^2(\Omega^a),}\\\\
\displaystyle{
\frac{h_{n_i}^{\frac{1}{2}}}{r_{n_i}^{\frac{N-1}{2}}}
u_{n_i,k}^b\rightharpoonup u_k^b\hbox{ weakly in
}H^1(\Omega^b)\hbox{ and strongly in }L^2(\Omega^b),}
\end{array}\right. \quad\forall k \in \mathbb N,\end{equation}
\begin{equation}\label{limitedeboleautovettoribiszhu}\left\{\begin{array}{l}\displaystyle{
\frac{1}{r_{n_i}}D_{x'} u_{n_i,k}^a\rightharpoonup \xi^a_k\hbox{
weakly in
}(L^2(\Omega^a))^{N-1},}\\\\
\displaystyle{
\frac{1}{h_{n_i}}\frac{h_{n_i}^{\frac{1}{2}}}{r_{n_i}^{\frac{N-1}{2}}}\partial_{x_N}
u_{n_i,k}^b\rightharpoonup \xi^b_k\hbox{ weakly in
}L^2(\Omega^b),}
\end{array}\right. \quad\forall k \in \mathbb N,\end{equation}
as $i\rightarrow +\infty$.

To obtain $u^k\in V_\infty$, for every $k\in \mathbb{N}$, it
remains to prove that
\begin{equation}\label{uak=0}u_{k}^a(0)=0, \quad\forall k\in \mathbb{N}, \quad\hbox{ if }N=2.\end{equation}

By virtue of the first line in
(\ref{limitedeboleautovettori+infty}) and of the fact that
$u_{n_i,k}^a(x',0)= u_{n_i,k}^b(r_{n_i}x',0)$ for $x'$ a.e. in
$\omega$, for obtaining (\ref{uak=0}), it is enough to prove that
\begin{equation}\label{traccia=0}
\lim_i\int_\omega u_{n_i,k}^b(r_{n_i}x',0)dx'=0,\quad\forall k\in
\mathbb{N}.
\end{equation}

By taking into account that $u_{n,k}^b\rightarrow0$ strongly in
$H^1(\Omega^b)$ and that $\displaystyle{\left\Vert\frac{\partial
u_{n,k}^b}{\partial x_N}\right\Vert_{L^2(\Omega^b)}\leq
c{r_{n}^{\frac{N-1}{2}}}h_n^{\frac{1}{2}}}$, by adapting the
 the proof of
(2.4) in \cite{GGLM1}, limit (\ref{traccia=0}) can be achieved, if
$N=2$.

 By choosing in (\ref{eigenproblemperun,k+infty}), written with
 $n=n_i$,
$$v=\left\{\begin{array}{ll}0& \hbox{ in }\omega\times
]\varepsilon_i,1[,\\\\\displaystyle{\frac{r_{n_i}^{\frac{N-1}{2}}}{h_{n_i}^{\frac{1}{2}}}v^b(r_{n_i}x')\frac{\varepsilon_i-x_N}{\varepsilon_i}}&
\hbox{ in }\omega\times
]0,\varepsilon_i[,\\\\\displaystyle{\frac{r_{n_i}^{\frac{N-1}{2}}}{h_{n_i}^{\frac{1}{2}}}v^b}
&\hbox{ in }\Omega^b,\end{array}\right.$$ with $v^b\in
C^\infty_0(\omega)$ and $\{\varepsilon_i\}_{i\in\mathbb N }\subset
]0,1[$, it results that
\begin{equation}\nonumber\begin{array}{l}\displaystyle{\frac{r_{n_i}^{\frac{N-1}{2}}}{h_{n_i}^{\frac{1}{2}}}
\int_{\omega}\int_0^{\varepsilon_i}\frac{1}{r_{n_i}}D_{x'}
u^a_{{n_i},k}(D_{x'}v^b)(r_{n_i}x')\frac{\varepsilon_i-x_N}{\varepsilon_i}-\partial_{x_N}u^a_{{n_i},k}v^b(r_{n_i}x')\frac{1}{\varepsilon_i}dx+\int_{\Omega^b}\frac{h_{n_i}^{\frac{1}{2}}}{r_{n_i}^{\frac{N-1}{2}}}D_{x'}
u^b_{{n_i},k}D_{x'}v^b
dx}\\\\
\displaystyle{=\lambda_{{n_i},k}\frac{r_{n_i}^{\frac{N-1}{2}}}{h_{n_i}^{\frac{1}{2}}}
\int_{\omega}\int_0^{\varepsilon_i}u^a_{{n_i},k}v^b(r_{n_i}x')\frac{\varepsilon_i-x_N}{\varepsilon_i}dx+\lambda_{{n_i},k}\int_{\Omega^b}\frac{h_{n_i}^{\frac{1}{2}}}{r_{n_i}^{\frac{N-1}{2}}}
u^b_{{n_i},k} v^bdx},\,\,\forall i\in \mathbb N, \quad\forall k\in
\mathbb N.
\end{array}\end{equation}
Consequently, by passing to the limit as $i\rightarrow +\infty$,
and by making use of (\ref{q=+infty}) and
(\ref{convlambda+infty})$\div$(\ref{limitedeboleautovettori+infty}),
one achieves
\begin{equation}\label{equainOb+inftyaou}
\int_{\omega}D_{x'}u^b_kD_{x'}v^bdx'=\lambda_k\int_{\omega}u^b_kv^bdx',\quad\forall
v^b\in H^1_0(\omega),\quad\forall k \in \mathbb N,
\end{equation}
if  the sequence $\{\varepsilon_i\}_{i\in\mathbb N }$ has been
chosen such that $\displaystyle{\lim_i\varepsilon_i= 0}$  and
$\displaystyle{\frac{r_{n_i}^{{N-1}}}{h_{n_i}}<< \varepsilon_i}$.

Assume now $N=2$.

By passing to the limit, as $n_i\rightarrow +\infty$, in
(\ref{eigenproblemperun,k+infty}), with
$v=\left\{\begin{array}{ll}v^a& \hbox{ in }\Omega^a,\\0 &\hbox{ in
}\Omega^b\end{array}\right.$ and $v^a\in H^1_0(]0,1[)$, and by
making use of (\ref{convlambda+infty}) and of the first line of
(\ref{limitedeboleautovettori+infty}), it turns out that
\begin{equation}\label{equainOa+infty}\vert\omega\vert
\int_0^1\partial_{x_N}u^a_k\partial_{x_N}v^adx_N=\lambda_k\vert\omega\vert\int_0^1u^a_kv^adx_N,\quad\forall
v^a\in H^1_0(]0,1[), \quad\forall k \in \mathbb N.
\end{equation}

 By adding
(\ref{equainOb+inftyaou}) to (\ref{equainOa+infty}), one obtains
that
\begin{equation}\label{problemalimite+inftyas}\left\{\begin{array}{l}
\displaystyle{u_{k}\in { V_\infty},}\\\\
\displaystyle{\alpha_1(u_{k}, v)=\lambda_{k}[u_{k}, v]_1,\quad
\forall v\in { V_\infty},}
\end{array}\right.\quad\forall k\in \mathbb N.\end{equation}
if $N=2$. Finally, by passing to the limit, as $n_i\rightarrow
+\infty$, in (\ref{deltahk+infty}), and by making use of
(\ref{limitedeboleautovettori+infty}) and
(\ref{problemalimite+inftyas}), it turns out that
\begin{equation}\label{3q1}
\alpha_1(\lambda_{k}^{-\frac{1}{2}}u_{k},\lambda_{h}^{-\frac{1}{2}}
u_{h})=[u_{k}, u_{h}]_1=\delta_{h,k},\quad\forall k,h\in \mathbb
N.\end{equation} Moreover,  by arguing as in  the proof of Theorem
\ref{maintheorem}, one proves that
\begin{equation}\label{ag}\displaystyle{\lim_{k}\lambda_{k}=+\infty.}\end{equation}

Now, assume $N\geq 4$.

Let $\{\varphi_n\}_{n\in \mathbb{N}}\subset C(\mathbb{R}^{N-1})$
be a sequence such that, for every $ n\in \mathbb{N}$,
 $\varphi_n=1$ in
$\{x'\in \mathbb{R}^{N-1}: \vert x'\vert\leq r_n\}$, $\varphi_n=0$
in $\mathbb{R}^{N-1}-\{x'\in \mathbb{R}^{N-1}: \vert x'\vert<
2r_n\}$, $\{\varphi_n\}_{n\in \mathbb{N}}\subset C^1(\{x'\in
\mathbb{R}^{N-1}: r_n\leq \vert x'\vert\leq 2r_n\})$, $0\leq
\varphi_n\leq 1$, $\displaystyle{\int_{\{x'\in \mathbb{R}^{N-1}:
r_n<\vert x'\vert< 2r_n\}}\vert D \varphi_n\vert^2dx'}= c
r_n^{N-3}$, where $c$ is a constant independent of $n$, but
dependent on $N-1$ (for the existence of such sequence, see
 Proposition  3.1 in \cite{GGLM1}).

 Now, by passing to the limit, as $n_i\rightarrow +\infty$, in
(\ref{eigenproblemperun,k+infty}), with
$v=\left\{\begin{array}{ll}v^a& \hbox{ in
}\Omega^a,\\v^a(0')\varphi_n(\frac{x'}{d}) &\hbox{ in
}\Omega^b\end{array}\right.$ and $v^a\in H^1(]0,1[)$, $v^a(1)=0$,
where $d=dist(0',\partial\omega)$, and by taking into account the
additional assumption (\ref{ulterioreipotesi}), it is easily seen
that
\begin{equation}\label{equainOa++infty}\vert\omega\vert
\int_0^1\partial_{x_N}u^a_k\partial_{x_N}v^adx_N=\vert\omega\vert\lambda_k\int_0^1u^a_kv^adx_N,\quad\forall
v^a\in H^1(]0,1[), \quad\forall k \in \mathbb N.
\end{equation}

By adding (\ref{equainOa++infty}) to (\ref{equainOb+inftyaou}),
one obtains (\ref{problemalimite+inftyas}) (and consequently
(\ref{3q1}) and (\ref{ag})), if  $4\leq N$ and $r_n^{N-1}\ll
h_n\ll r_n^2$.

 If $N=3$, one obtains (\ref{equainOa++infty}) by arguing as
 above, but by choosing $\varphi_n$ such that, for every $ n\in \mathbb{N}$,
 $\varphi_n=1$ in
$\{x'\in \mathbb{R}^{N-1}: \vert x'\vert\leq r_n\}$, $\varphi_n=0$
in $\mathbb{R}^{N-1}-\{x'\in \mathbb{R}^{N-1}: \vert x'\vert<
\sqrt{r_n}\}$, $\{\varphi_n\}_{n\in \mathbb{N}}\subset C^1(\{x'\in
\mathbb{R}^{N-1}: r_n\leq \vert x'\vert\leq \sqrt{r_n}\})$, $0\leq
\varphi_n\leq 1$, $\displaystyle{\int_{\{x'\in \mathbb{R}^{N-1}:
r_n<\vert x'\vert< \sqrt{r_n}\}}\vert D \varphi_n\vert^2dx'}= -c
(\log r_n)^{-1} $, where c is a positive constant (for the
existence of such sequence, see
 Proposition  3.1 in \cite{GGLM1}), and  $r^2 \ll h_n\ll -r_n^2\log r_n$.

 The identification of $u^a_k$ remains an open question when
 $N\geq 3$ and assumption (\ref{ulterioreipotesi}) is not satisfied.

Let, now, $N=2$, or  $N=3$ and $r_n^{2}\ll h_n\ll -r_n^2\log r_n$,
or  $ N\geq 4$ and $r_n^{N-1}\ll h_n\ll r_n^2$. As in Step 2 of
the proof of Theorem \ref{maintheorem}, by arguing by
contradiction, one can show that   there not exist $(\overline u,
\overline \lambda)\in V_\infty\times\mathbb{R}$ satisfying the
following problem:
\begin{equation}\label{probleoverlinelambdakj1}\left\{\begin{array}{l}
\displaystyle{\overline
u\in { V_\infty},}\\\\
\displaystyle{\alpha_1(\overline u, v)=\overline\lambda [\overline
u, v]_1,\quad
\forall v\in { V_\infty},}\\\\

[\overline u, u_k]_1=0,\quad\forall k\in\mathbb N\\\\

[\overline u,\overline u ]_1=1.
\end{array}\right.\end{equation}
Precisely, assume that there exists $(\overline u, \overline
\lambda)\in V_\infty\times\mathbb{R}$ satisfying
(\ref{probleoverlinelambdakj1}). Let
 ${\overline k}\in \mathbb{N}$ be such that
\begin{equation}\nonumber
\overline\lambda< \lambda_{\overline k},
\end{equation} and,
for every $n\in \mathbb{N}$, let $\varphi_n$ be the solution of
the following problem:
\begin{equation}\nonumber\left\{\begin{array}{l}
\displaystyle{\varphi_n\in { V}_n,}\\\\
\displaystyle{a_n(\varphi_n,
v)=\overline\lambda\left(\left(\overline
u^a,\frac{r_n^{\frac{N-1}{2}}}{h_n^{\frac{1}{2}}}\overline
u^b\right), v\right)_n,\quad \forall v\in { V}_n,}
\end{array}\right.\end{equation}
where $\overline{u}=(\overline u^a,\overline u^b)$. Then, it is
easy to prove that
\begin{equation}\nonumber\displaystyle{
\varphi_{n}^a\rightharpoonup \overline {u}^a\hbox{ weakly in
}H^1(\Omega^a),\quad
\frac{h_{n}^{\frac{1}{2}}}{r_{n}^{\frac{N-1}{2}}}\varphi_{n}^b\rightharpoonup
\overline{u}^b\hbox{ weakly in }H^1(\Omega^b),}
\end{equation}
as $n\rightarrow +\infty$. Then, by proceeding as in the proof of
Step 2 of Theorem \ref{maintheorem},  one reaches a contradiction.

To complete the proof, one can argue as  Step 3 of the proof of
Theorem \ref{maintheorem}.
\end{proof}

\textit{Proof of Theorem \ref{TH3}.}  Theorem \ref{TH3} follows
immediately from Theorem \ref{theorem+infty}, by arguing as in the
proof of Theorem \ref{TH1} at the end of Section \ref{mainresult}.
$\hfill\Box$

\section*{Acknowledgments}
The authors wish to thank  "LATP" of the University of Provence
(Aix-Marseille I), where the first author was invited as visiting
professor in May 2005 and in May 2006 and where this research has
been initiated and concluded.

This paper is part of the project: "Strutture sottili" of the
program 2004-2006: "Collaborazioni interuniversitarie
internazionali" of the Italian Ministry of Education, University
and Research.



\end{document}